\documentclass[aps,prc,nofootinbib,showpacs,superscriptaddress,groupedaddress]{revtex4-2}
\usepackage{subfiles}
\usepackage{amsmath,amssymb,amsbsy,bm}
\usepackage{graphicx}
\usepackage{comment}
\usepackage{float}
\usepackage[colorlinks=true,linkcolor=blue,citecolor=blue,urlcolor=blue]{hyperref}

\newcommand{\unity}{\mathbb{I}}
\newcommand{\E}{\mathrm{E}}
\newcommand{\Var}{\mathrm{Var}}
\newcommand{\Cov}{\mathrm{Cov}}

\begin{document}

\title{Optimally Emulating Smooth Functions for Bayesian Analyses}
\author{Scott Pratt, Oleh Savchuk, Eren Erdogan, Ekaksh Kataria}
\affiliation{Department of Physics and Astronomy and National Superconducting Cyclotron Laboratory\\
Michigan State University, East Lansing, MI 48824~~USA}
\author{Christal Martin}
\affiliation{Department of Physics, University of Tennessee, Knoxville, TN 37996}
\author{Syed Afrid Jahan}
\affiliation{Department of Physics and Astronomy, Wayne State University, Detroit, MI 48202}
\author{Frederi Viens}
\affiliation{Department of Statistics, Rice University, Houston, TX  77005}

\date{\today}


\begin{abstract}
For Bayesian analyses aimed at constraining the parameter space of computationally expensive models by comparison with heterogeneous data sets, model emulators play a pivotal role. Model emulators enable searches through high-dimensional parameter spaces, searches which might otherwise require millions of full model runs. Here, an emulation strategy is presented based on the assumption that the model output can be well described by a polynomial expansion with all terms being independent of one another and with higher order terms being given decreasing weights according to a convergence parameter. It is found that if the expansion's form is constrained to being invariant to rotation and translation of the parameter space that the polynomial form becomes equivalent to that of a Gaussian process emulator using a correlation kernel with a Gaussian form. In the limit that the function is smooth, quantified by a convergence parameter, $\Lambda$, it is found that rather few training points are required to reproduce the full model and provide a good estimate of its uncertainty.  An accuracy metric is defined based on an average of the expected uncertainty over the prior before training data are actually gathered. Analytic expressions for the expected accuracy are provided based solely on an expected value of $\Lambda$ and the location of the training points. This enables a tenable optimized placement of training points, which is shown to be significantly superior to random placement or that from Latin hyper-cube sampling. Examples of optimized placement are provided, and the accuracy is studied as a function of various variables, such as the number of model parameters, the number of training points, and the convergence parameter.  
\end{abstract}

\maketitle





\setcounter{section}{0}
\section{Introduction}

In the quest to provide rigorous quantitative constraints of high-dimension model-parameter spaces by comparing sophisticated computationally expensive numerical models to experimental results, one must sample the model-parameter space by running the full-model at large numbers of points. For example, sampling a 12-parameter model with only three points in each direction requires $3^{12}\approx 500,000$ runs. Thus, Markov Chain Monte Carlo (MCMC) searches, weighted by the likelihood that the model reproduces the experiment, often involve millions of full-model evaluations. If a single run of the full-model at one point in parameter space requires days of CPU time, the analysis becomes untenable. Model emulation provides a solution to this problem. If one can train an emulator to reproduce the full-model by training it on only a few hundred full-model runs or fewer, and if the emulator can provide a reasonable estimate of its uncertainty, the emulator can take the place of the full model in the MCMC search, and thus enable the analysis. 

\subsection{Motivation, context, and related work}\label{sect:context}

This paper is motivated by the theory and practice of heavy-ion physics, a field where emulators have been applied to various large-scale analyses, often involving more than a dozen model parameters, with models that often require days of CPU time. Contextual references in heavy-ion physics abound: 
\cite{Domingues:2024pom,JETSCAPE:2024cqe,JETSCAPE:2024ofm,Paquet:2023rfd,JETSCAPE:2023ikg,Heffernan:2023gye,Heffernan:2023utr,Weiss:2023yoj,Ji:2022xzo,Heffernan:2022swr,Liyanage:2022byj,JETSCAPE:2021rfx,JETSCAPE:2021ehl,Bernhard:2019bmu,Moreland:2018gsh,Moreland:2018jos,Xu:2017obm,Ke:2017xyi,Xu:2017nmc,Auvinen:2017fjw,Xu:2017hgt,Auvinen:2017pny,Bernhard:2017vql,Ke:2016jrd,Auvinen:2016tgi,Bernhard:2015hxa,Jahan:2024wpj,JETSCAPE:2023nuf,Shen:2023pgb,Shen:2023pgb,Mantysaari:2022ypp,Mantysaari:2022ffw,JETSCAPE:2020mzn,JETSCAPE:2020avt,Pratt:2015zsa,Sangaline:2015isa,Novak:2013bqa,Li:2025iqq,Alqahtani:2024ejg,Wang:2024ahj,Li:2023ydi,Giannini:2022bkn,Jia:2013tja,Morfouace:2019jky}.
For example, in a global analysis of 14 model parameters \cite{Pratt:2015zsa}, Gaussian-processes (GP)  \cite{RasmussenGP,WilliamsRasmussen96} were used as emulators of the output from hybrid hydrodynamic/Boltzmann approaches with accuracies near one percent, after being trained on  1200 points. One might ask whether one can do as well, or better, with fewer training points. The field of heavy-ion physics is similar to many other physical or biological fields of application, in that  the relation between the model parameters and model  output tends to be roughly linear in a physical parameter range. For example, in high-energy Lead-Lead (Pb-Pb) collisions at the Large Hadron Collider (LHC), the mean momentum of outgoing protons rises monotonically with the stiffness of the equation of state. That dependence should not be perfectly linear, but one would expect that, within a reasonable prior parameter range, the dependence is well described by linear or quadratic functions, and that higher-order terms should become increasingly negligible. This expectation  also applies to most any other observable with respect to nearly any other model parameter. For lower-energy collisions, where phase transitions are perhaps more relevant, one then imagines a dependence that has significant structure, e.g. a clear threshold w.r.t. some parameter or set of parameter hypersurfaces, with roughly linear behaviors in the regions determined by these surfaces. At the higher LHC energies, these complications would not occur, and the methods in this paper would apply directly to the entire physical parameter range. 

As a quick sanity check, putting aside  complete UQ considerations, imagining a scalar model, with $N_{\rm pars}$ parameters, if it were nearly perfectly linear, a first-order polynomial (linear approximation) would provide an excellent fit, and would presumably result in acceptable UQ bands. One could then run the full model only $N_{\rm pars}+1$ times, at points fairly regularly far apart, to constrain the linear approximation's $N_{\rm pars}$ slopes and the intercept. For a 12-parameter model, this would  require only 13 full model runs. Imagining  that the model was instead nearly perfectly quadratic, one would need only $N_{\rm pars}(N_{\rm pars}+1)/2$ additional full-model runs to constrain all the quadratic contributions. For the 12-parameter case, this would lead to 91 full model runs, with an expectation that this way of contraining the model would still lead to acceptable UQ bands. This is fewer than the number of full-model runs often used in recent analyses in heavy-ion physics, by orders of magnitude, for modeling scenarios where there is arguably little benefit, within physical parameter ranges, of a cubic approximation over a quadratic one.

Another way to say this, which is prudent and parsimonious, both in computational and modeling terms, is that in any domain of application, such as heavy-ion physics, where most models feature rather smooth responses to their parameters, is that one should reconsider polynomial fits, and investigate the degree to which they might, or might not, improve upon the GP approaches cited above. The first goal of this paper is to exploit the expectation that a polynomial fit should be adequate  systematically, by presenting a way to construct a polynomial expansion where the importance of various terms in the expansion decreases with increasing order. In the physics literature, this expectation of adequacy (of a polynomial approximation) is often labeled as ``convergence'' of the approximation, where the asymptotic notion in this convergence is with respect to the degree  of the polynomial as it increases to infinity, even though practitioners talk about convergence even if the degree is small, as long as the error in the approximation fulfills the practitioners' ``expectations'', which in this case should be understood as their stated or implicit tolerances.  To be able to declare that approximations are ``converging'' in the sense that the errors are tolerable, it is best to quantify those expectations in the most transparent terms possible, which is often done in terms of a ``convergence'' parameter, $\Lambda$. Once one has agreed on such a quantitative framework,  one can work to define and find an optimum set of expansion coefficients, along with a quantification of their uncertainty. Herein, the emulator is built on a polynomial expansion where the polynomial's variables are the $N_{\rm pars}$ values $\theta_i/\Lambda$, where the model parameters are $\theta_i;i=1,\cdots,N_{\rm pars}$, and are all scaled so that they vary of order unity within the prior. Thus, higher-order terms contribute as negative integer powers of the convergence parameter $\Lambda$.

Emulators defined by choosing polynomial fits based on least-square minimization algorithms have long been viewed as  good classes of candidates. More modern approaches improve on such fits by providing estimates that exactly reproduce the training values when evaluated at the training points, and provide an estimate of the emulator's uncertainty as a function of the vector of parameters $\theta$. These approaches which reproduce training values are typically based on GP assumptions, where the predictions are all based on an assumed form of the covariance between values $Y$ and $Y'$ of the model evaluated at two parameter-vector points $\theta$ and $\theta'$. That correlation kernel might be any function of $\theta$ and $\theta'$, but a particularly simple and popular form is a ``Gaussian'' covariance kernel, $\Cov_0(Y_1,Y_2|\sigma_A,\Lambda)=\sigma_A^2e^{-|\theta - \theta'|^2/(2\Lambda^2)}$. Here, the parameters $\sigma_A$ and $\Lambda$ are referred to as hyper-parameters, to differentiate them from model parameters.  They are used only to tune the emulator, and their tuning can be optimized based on training data. Here, the subscript $0$ denotes the fact that the covariance is taken with respect to the prior distribution. It has been shown that any polynomial expansion, or any expansion over some set of basis functions, where the expansion's terms have random coefficients, can be optimized to reproduce training points. These can then be equated to GP forms once one defines the appropriate correlation kernel \cite{RasmussenGP}. Previous approaches have typically been constructed from sets of orthogonal basis functions, e.g. \cite{Mak24}, but this requirement is not necessary and is not enforced herein. The sum over polynomial terms with random coefficients is sometimes referred to as ``polynomial chaos'' \cite{Xiu2010}. Once trained, one can optimize all the random coefficients for the polynomial expansion. One can then calculate the emulator values at any parameter vector $\theta$ using those optimal coefficients. However, there are many more coefficients than training constraints, since one has the option to truncate the full Taylor series at an arbitrarily high degree. Fortunately, once one has written an equivalent form for the correlations,  the emulator can much more efficiently be expressed in terms of the correlation kernel and one can avoid performing sums over all monomial terms. This is sometimes referred to as the {\it kernel trick} \cite{RasmussenGP}.

\subsection{Narrative summary of this paper's contributions}

In this work's first part, we consider emulators as polynomial expansions with random coefficients, and we constrain the emulators' form by requiring that the emulating behavior be invariant by rotation and by translation. Under this requirement, we find that the equivalent correlation kernel must then have a simple Gaussian form. Thus, if one believes that the model is equally sensitive to all parameters and that there is no reason to expand the polynomials about one point rather than around some other point, then the equivalent GP correlation kernel must be Gaussian, and described via the form involving the two hyper-parameters $\sigma_A$ and $\Lambda$ mentioned above. The polynomial description then becomes equivalent to a GP description with a correlation kernel with a simple Gaussian form. If the user believes that the sensitivity w.r.t. some model parameter should be less than that of another parameter, rather than applying different convergence parameters, one can scale the parameters so that the expected sensitivity is isotropic in the scaled parameter basis. Thus, the assumption of the model having a rotation-invariant behavior can reasonably apply to many problems after appropriate rescaling. The assumption of translation invariance is one of spatial stationarity across parameter space, which is also a strong modeling preference for most problems. In particular, it significantly reduces the statistically grounded UQ on any parameter, avoiding the risk of overfitting a model for UQ itself. 

On the strength of the translation and rotation invariance assumption of the model's parameter-space sensitivity, the single scalar convergence parameter, $\Lambda$,  sets the relative mutual importance of terms of different orders in the polynomial expansion. We find that $\Lambda$ is  precisely the same as the correlation length in the GP interpretation with a Gaussian correlation kernel. For the polynomial expansion to ``converge'', i.e. for its truncation to yield a tolerable error, within a few orders, at a minimum, one must have the convergence parameter $\Lambda$ larger than the characteristic scale of the prior on parameter space, and this would typically not be enough for the truncation error to be tolerable. If the model parameters are scaled to vary roughly between $-1$ and $+1$, then $\Lambda$ should be significantly larger than unity if the expansion is to be well defined by a few terms, or equivalently, to be constrained via calibration by a few training points. After translating to the equivalent GP paradigm, this is equivalently stated as having the GP correlation length, which coincides with $\Lambda$, being significantly larger than the size of the prior.

Conversely, if one assumes a Gaussian correlation kernel in a GP context, the relations derived here identify the form of the corresponding polynomial expansion. Thus, one can use this form to generate random functions that are consistent with the Gaussian kernel, albeit with some cutoff in the order of the expansion's truncation. This should be useful for testing emulation software, as it should allow one to build models that correspond to the desired correlation kernel and should, by construction, have uncertainties that are well understood and quantified. 

Our work's second thrust is to understand, for expansions which are rapidly converging (where only a small number of terms are needed to get tolerable truncation errors), how one should optimally choose the points at which the full model is run to train the emulator. We present a metric that represents the expected accuracy (expected truncation error) of the emulator as a function of the convergence parameter, $\Lambda$, and as a function of the positions of the training points. For a fixed $\Lambda$, optimized arrangements of training points are then generated for various cases involving either Gaussian or uniform priors on the parameters. For test cases, this procedure is found to be significantly superior to either random placement of the training points throughout the prior, or to Latin hyper-cube sampling \cite{lhs}. In this setting of rapidly converging expansions, one might expect that the number of training points might be chosen to match that needed for linear, quadratic or cubic fits. For a 6-parameter example, that would require 7, 28 or 84 full-model runs respectively. Here, we show that the marginal improvement to the emulation decreases at each of the two corresponding thresholds, assuming that the training points are optimally chosen. In practice, this analysis would provide a user who has a clear understanding of their computational costs, and who has a clear quantification of their aversion to error risk, which can simply take the form of a worst-case scenario cost of tolerance, to run an elementary cost-benefit analysis.  

By accomplishing the goals above, one should be able to gain a better understanding of how many model runs might be necessary to design an emulator to the desired accuracy, where the term ``desired'' can be built into a cost-benefit analysis of the aforementioned ilk. This decision framework will still depend on some pre-estimate or other choice of the convergence parameter $\Lambda$, but it turns out, from our testing, that the actual placement of the training points appears rather insensitive to that estimate for $\Lambda$, even though the accuracy is strongly sensitive to the estimate. For users who knows ahead of time that some of the model parameters are more likely to drive the functional behavior than others, we provide a seamless method to exploit that prior knowledge. This can assist with overcoming the challenges of high dimensionality. We find herein that such flexibility turns out to be critical in determining how confidently one might ultimately emulate models.

In the next section, we present the mathematical form of such an emulator, and a derivation of the optimum set of coefficients, given the training data. Closed expressions for the emulator uncertainty are also derived.  The emulator's form also accounts for  models with noisy training, i.e. those for which point-by-point uncertainty is acknowledged and taken into account, i.e. where that uhcertainty is propagated in a fashion consistent with the modeling. The physical arguments for choosing the specific emulator form are discussed. The equivalence with a GP emulator is established carefully. The subsequent section presents a tractable method to optimize the location of training points in the $N_{\rm pars}-$dimensional parameter space. For the case of Gaussian priors with $N_{\rm pars}+1$ training points, the placement follows a straightforward $N_{\rm pars}-$dimensional simplex, but for most cases the optimized solution does not have a simple symmetry. Instead, our method finds the optimal points numerically. In Section \ref{sec:testing} we test our emulator procedure by applying it to some models with random polynomial expansions or randomized trigonometric or exponential forms. This allows us to gain a solid quantified grasp for how well one might, in practice,  not only emulate a function, but also understand how well one might trust the fidelity of our emulator's quantitatively estimated uncertainty. This faithfulness is studied as a function of the dimensionality of the model's parameter space and of the convergence parameter $\Lambda$.

\setcounter{section}{1}
\section{Theoretical Foundation of Emulating Convergent Functions}\label{sec:theory}

The form and basic properties of our proposed class of emulators is described in detail in this section. The form is that of a polynomial expansion in terms of the model parameters, i.e. expanded as a multivariate power series with respect to the parameter vector, where the coefficient for each term in the polynomial is given an independent centered Gaussian (mean-zero normal) weight. Since centered Gaussian variables have a scale property with respect to their standard deviation, it is convenient and elementary to standardize the Gaussian coefficients by pulling out their scales as explicit deterministic emulation weights.  Using these weights, our form, while general, is easily constrained to assume that terms of higher order should contribute less. For example, in the absence of training data, quadratic terms should be expected to contribute less than linear terms, and cubic terms should contribute less than quadratic terms. What separates the form presented here compared to previous studies is that we also constrain the expansion to behave invariantly by rotation and by translation in parameter space. Thus, moving the origin about which the polynomial expansion is constructed, or rotating the coordinate system about any such origin, should not affect the emulator. 

In fact, in the next subsection we show that this natural constraint implies a precise and unique form for the deterministic emulator weights, well illustrating the intuition that weights should decrease as terms' order increases. We shall see that the form of the resulting polynomial expansion depends only on two emulator parameters, $\sigma_A$, which sets the strength of the fluctuation of the pre-trained emulator at any point, and $\Lambda$, which drives the convergence of the polynomial expansion, as it is a power series expansion in terms of powers of the scaled parameter vector $\Lambda^{-1}\vec{\theta}$. A reader might pause at this point to recognize that we will be dealing with $N_\text{pars}$ physical model parameters assembled in the vector $\vec{\theta}$ and with two scalar emulator parameters $\sigma_A$ and $\Lambda$.

It is well known that expansions, including polynomial expansions, which are linear in their random coefficients, lead to emulators whose behavior can easily be expressed in terms of the pre-trained covariance kernel, $C_0(\vec{\theta}_1,\vec{\theta}_2)$. Both the emulator and its error are easily expressed in terms of $C_0$, which allows one to avoid calculations where one must sum over the actual terms of the polynomial expansion. This result is known as the {\it kernel trick} \cite{RasmussenGP}. As we said, this is not new. The new constructive element we introduce here, is that our assumption that all the random terms in the polynomial expansion are independent of each other, and our constraint of rotational and translational invariance, together imply a unique form for the expansion, and thus a unique form for the kernel. That unique form turns out to be that of an ordinary Gaussian kernel (i.e. a simple quadratic exponential decay) with $\Lambda$ playing the role of characteristic length, in parameter space, of that quadratic exponential correlation-decay. By standardizing all parameters to be of the same magnitude, one then only requires a single characteristic length, i.e. just one scalar $\Lambda$.

Throughout this section the model parameters, $\vec{\theta}$, are assumed to have been translated, rotated and scaled so that their expected sensitivities to a given step size $\Delta\Theta$ are similar. For example, if a model parameter $x_i$, has an initial prior, $a_i<x_i<b_i$, one defines a new translated and scaled variable
\begin{eqnarray}
\theta_i=\frac{2\theta_{{\rm max},i}}{(b_i-a_i)}(x_i-(a_i+b_i)/2).
\end{eqnarray}
where now the prior on $\theta_i$ is uniform in the symmetric range $\pm\theta_{{\rm max},i}$. The flexibility with choosing the scale $\theta_{{\rm max},i}$ is use so that, for each $i$, one would not expect the observable, $y(\vec{\theta})$, to be any more sensitive to steps of size $\Delta\Theta$ in the direction of $\theta_i$ than to steps of the same size in the direction of some other parameter $\theta_j$. Thus, if the sensitivity to $\theta_1$ was expected to be higher than the sensitivity to $\theta_2$, one would reduce the ratio of $\theta_{{\rm max},2}$ to $\theta_{{\rm max},1}$ until the expected sensitivities were similar. All this can be attempted mathematically when partial derivatives of the model for $y(\vec{\theta})$ can be estimated with some confidence, or it can be done empirically thanks to a limited number of full-model numerical evaluation, not emulation. Further, on the question of fixing issues with rotational dependence, to consider a simple example, if one expected  $y(\vec{\theta})$ to be more sensitive to a given direction in parameter space, e.g., $\theta_1+\theta_2$, compared to some orthogonal direction, $\theta_1-\theta_2$, one could rotate the coordinate system by 45 degrees, then scale the two dimensions accordingly. More generally, one can define a matrix, $S_{ij}$, describing the expected or empirical variabiity of $y$ with respect to $\Delta\theta_i \Delta\theta_j$, then redefine a coordinate system where the new parameters were eigenvector of $S$, then scale the new parameters so that the eigenvalues of $S$ were similar. One would then perform the model emulation in this new basis where the functions was not expected to vary more in any specific direction for the same step size, although the range of the prior could vary per directions. For the case where one had no information indicating that the model were more likely to vary across the prior in one direction than in another, all the values of $\theta_{{\rm max},i}$ would all be set to equal $1$. 

As we mentioned, when using a uniform prior on each parameter, we shall assume that the parameter $\theta_i$ is distributed according to the following probability density:\begin{equation}\label{eq:uniformpriorsdef}
P_i(\theta_i)=\left\{
\begin{array}{rl}
\frac{1}{2\theta_{{\rm max},i}},&-\theta_{{\rm max},i}<\theta_i<\theta_{{\rm max},i},\\
0,&{\rm otherwise}
\end{array}
\right.
\end{equation}
and as we just said, if all parameters are considered to be of equal importance, we let $\theta_{{\rm max},i}=1$ for every $i$. If the model function $y(\vec{\theta})$ is expected to vary mostly along a specific $i$th parameter's direction, since everthing is invariant to a single scaling factor, we need only scale that parameter with $\theta_{{\rm max},i}=1$, and scale other parameters with appropriately chosen smaller values of $\theta_{{\rm max},j},j\neq i$. 

If the parameters have Gaussian priors, one can scale and translate each of them so that the new marginal priors are Gaussian of the form
\begin{eqnarray}\label{eq:gausspriorsdef}
P_i(\theta_i)&=&\frac{1}{(2\pi/\beta_i)^{1/2}}e^{-\beta_i\theta_i^2/2}.
\end{eqnarray}
where the $\beta_i$'s are chosen to place the parameter sensitivities on presumed equal footing. But the Gaussian case is easier to treat formally than the uniform case, because a multivariate Gaussian distribution has an intrinsic property of rotational invariance that the multivariate uniform distribution lacks. Therefore, if one insists on a specific multivariate Gaussian prior distribution on the entire vector $\vec{\theta}$, with mean $\vec{\mu}$ and covariance matrix $S$, one then simply shifts the entire vector $\vec{\theta}$ by subtracting $\mu$, and one rotates and scales it simultaneously by multiplying it by a Cholesky square root $L$ of $S$, i.e. a lower-triangular matrix $L$ such that $S=LL^{\top}$. The new resulting multivariate prior is then the product of all the normal densities in \eqref{eq:gausspriorsdef} with each $\beta_i=1$. These $\beta_i$'s can then also be adjusted to make up for uneven sensitivites, similarly to what we suggested with the $\theta_{{\rm max},i}$ scales for uniform priors.

Again in the multivariate normal case, if we now believe that all the parameters are put on equal footing, but we want this multivariate prior to have the same dispersion properties as in the uniform case, one might want to match variances. The variance of the uniform distribution on $[-1,1]$ being equal to $1/3$, the value of all the $\beta_i$'s need only be set to equal $3$ as well for the variances to match in the Gaussian and uniform priors cases. As we mentioned in the uniform prior case, if some given parameter (\# $i$) is expected to drive most of the variability that parameter is assigned $\beta_i=3$ and all other parameters are assigned lower values for their $\beta_j$'s. Certainly, matching variances is not the only way to match dispersions. If instead one wanted a density in \eqref{eq:gausspriorsdef} with a $\beta_i$ which matches the support of the uniform distribution $[-1,1]$ at a confidence/credibility level of $95\%$, then one will want to choose not $\beta_i=3$ but rather closer to $\beta_i=4$ (strictly speaking, it would be closer to $\beta_i=3.84$.) This is because the normal distribution has higher moments that increase slightly faster than those of the uniform distribution. In the sequel, we will stick with the variance-matching choice of $\beta_i=3$ for our expository purposes.   

With our choice to use polynomial expansions in terms of $\theta_i/\Lambda$, as discussed above, it may be the case that the function to be emulated is known to be less sensitive to some model parameters than it is to others. Rather than assigning different $\Lambda$ for different model parameters to reflect $\theta_{{\rm max},i}$ and $\beta_i^{-1/2}$ difference in sensitivities, the prior scales $\theta_{{\rm max},i}$ and $\beta_i^{-1/2}$ are modified from the values of $\theta_{{\rm max},i}=1$ and $\beta_i^{-1/2}=1/\sqrt{3}$. This choice makes the mathematical description of the emulator simpler and maintains rotational invariance, while not forfeiting the ability to assume some parameters contribute in a more convergent manner than others. Because the emulator is a function of $\theta_i/\Lambda$, scaling $\theta_i$ is mathematically identical to scaling $\Lambda$. Moreover, one can assume different prior forms, e.g. half Gaussians, exponentials, Laplace distributions, and these choices do not affect any of the derivations or results in this section. However, the choice of the parameter priors are critical in Section \ref{sec:optimize}, where methods to optimize the training-point placement are considered.

The remainder of this section includes several subsections briefly reviewing the methods needed to tune and train our proposed emulator. Training is accomplished by calculating the full function $y(\vec{\theta})$ at $N_{\rm train}$ predetermined locations $\theta_a,a=1,2,\cdots ,N_{\rm train}$ in model-parameter space. The training data is thus comprised of knowing $y_a=y(\vec{\theta}_a)$ for $a=1,2,\cdots ,N_{\rm train}$. Expressions for the optimized coefficients are found in terms of the training data. After applying the aforementioned ``kernel trick'', generating these coefficients is mostly unnecessary given that the emulator and its uncertainty can readily be expressed in terms of the training data and the covariance kernel, without ever performing (inifinite) sums over all polynomial terms. Given that the resulting emulator is identical to one used in many GP treatments, our detailed mathematical and computational analysis of training and tuning might be viewed as unnecessary in practice. Nonetheless, it is insightful to understand what precise set of basis functions results in a GP emulator with a Gaussian covariance kernel. 

One part of our analysis of training and tuning which does not appear in many GP treatments is the inclusion of the emulator uncertainty due to the uncertainty in choosing the emulation hyper-parameter $\Lambda$. Although this uncertainty term is not usually significant once one has a reasonable distribution of training points, it plays an important role in estimating the number of training points required for a desired emulation accuracy, and in optimizing training-point placement itself, a subject considered in detail in Section \ref{sec:optimize}. To summarize briefly, one of the conclusions of our analysis of optimal training point placement is that standard methods like random or quasi-random placement, or latin-hypercube placement, are clearly sub-optimal for our emulator class, and therefore are suboptimal for the standard GP emulators with Gaussian kernels.

\subsection{General Form and the Assumption of Analyticity}

Recall that our $N_{\rm pars}$-dimensional model parameters can be expressed as a vector $\vec{\theta}$, and that we assume that the original model parameters have been rotated, shifted and scaled so that the prior ranges are defined as described above, with scale modulations to account for  higher presumed sensitivities in some directions than others in parameter space. The first assumption for choosing the functional form is that it be analytic, so that it can be expanded as a polynomial in full parameter space, at least in the range of the prior. To avoid abusing the complex-analysis term ``analytic'' here, we can say that its meaning is the classical one, for the power series to have an infinite radius of convergence when the prior's support is the full $N_{\rm pars}$-dimensional Euclidean space, and has a radius of convergence defining a sphere in that Euclidean space which contains the support of any prior when it is bounded.   Here, we assume the following  form for the emulator, $Y(\vec{\theta})$,
\begin{eqnarray}\label{eq:generalform0}
Y(\vec{\theta})&=&\sum_{\vec{n}} d_{\vec{n}}(|\vec{\theta}|^2) f_{K(\vec{n})}(|\vec{\theta}|^2)
A_{\vec{n}}
\left(\frac{\theta_1}{\Lambda}\right)^{n_1}\left(\frac{\theta_2}{\Lambda}\right)^{n_2}\cdots 
\left(\frac{\theta_N}{\Lambda}\right)^{n_{N{\rm pars}}}.
\end{eqnarray}
Throughout the paper, the uppercase function $Y(\vec{\theta})$ will refer to the emulator whereas the lowercase function $y(\vec{\theta})$ refers to the evaluation of the full model. Here, $K=K(\vec{n})=|\vec{n}|=n_1+n_2+\cdots n_{N_{\rm pars}}$ is the order of the polynomial, $d_{\vec{n}}$ and $f_{K(\vec{n})}$ are functions of the squared norm $|\vec{\theta}|^2$. The vector $\vec{n}$ refers to the $N_{\rm pars}$-dimensional vector comprised of non-negative integers, $(n_1,n_2,\cdots,n_k)$. 

The random coefficients $A_{\vec{n}}$ are considered to be independent in the absence of training data. Therefore, each polynomial term can be considered as a basis function. This set of basis functions is clearly non-orthogonal. The normalization and orthogonality properties are discussed in Sec. \ref{sec:basisprops}. The form of the emulator class $Y(\vec{\theta})$ appears rather general, in that all terms for a given $K$ are independent. For the emulator class $Y(\vec{\theta})$ to be an analytic function of emulator class $\vec{\theta}$, and to avoid any unnecessary contortions, the functions $f_K$ and $d_{\vec{n}}$ are best chosen to be polynomial expansions of the squared Euclidean norm $|\vec{\theta}|^2$. With this natural choice, the emulator $Y(\vec{\theta})$'s current form above, where the infinite series is assumed to be absolutely summable within an appropriate radius of convergence, spans all possible analytic functions of $Y(\vec{\theta})$, i.e. functions that can be written as polynomials of $\vec{\theta}$, or as uniform limits of these polynomials on compact sets in parameter space. However, this precludes explicit dependence on terms such as $\theta_1/|\vec{\theta}|$.

Our second assumption is about the random coefficients $A_{\vec{n}}$. We assume that the coefficients are independent and identically distributed, mean-zero Gaussian random variables, with common variance $\sigma^2_A$, i.e. with marginal density $\varphi_{A_{\vec{n}}}(a)$ and covariance structure given by
\begin{eqnarray}\label{eq:AGauss}
\varphi_{A_{\vec{n}}}(a)&=&\frac{1}{\sqrt{2\pi\sigma^2_A}} e^{-a^2/2\sigma_A^2},\\
\nonumber
\Cov_0(A_{\vec{m}},A_{\vec{n}})&=&\sigma^2_A\delta_{\vec{m},\vec{n}}.
\end{eqnarray}
We use the subscript ``0'' to refer to the averaging (via expectations or covariances) over the prior, i.e. without any constraints from matching training data. The particular choice of using Gaussian coefficients does not affect any of the derivations of this section unless one is considering moments of higher order than quadratic, e.g a correlation of order $Y^4$ such as $\E_0(Y(\vec{\theta}_a)Y(\vec{\theta}_b)Y(\vec{\theta}_c)Y(\vec{\theta}_d))$. In Sec. \ref{sec:Agauss} it will be shown that each independent random coefficient $A_{\vec{n}}$ must have a Gaussian distribution if the higher-order moments are to all be independent of translation and rotation. Any distribution that averages to zero and has some variance $\sigma_A^2$ would give the same behavior for moments involving only quadratic terms of $Y$. The choice that each random coefficient's distribution is characterized by the same single scale parameter $\sigma_A$ does not reduce the generality because $d_{\vec{n}}$ is still arbitrary up to this point. 

Similarly, contrary to what intuition might dictate, the choice that these coefficients are all independent of one another does not fundamentally reduce the generality of the expression except for the fact that we choose to work a unique random coefficient per monomial, disallowing repeats among the $\vec{n}$'s. To illustrate this point, assume for a moment that one allows such repeats. Then the above form does not eliminate the possibility that different terms could be correlated. Indeed, as long as the overall joint distribution of the (infinite-dimensional) vector of all coefficients is assumed to be that of a Gaussian vector, a sum of independent terms can be achieved as follows. Any Gaussian vector, even with countably infinite dimension, has a covariance matrix that can be expressed for instance as its Cholesky decomposition. Then, its corresponding square root can be inverted and multiplied by the Gaussian vector, resulting in an explicit series of independent Gaussian terms for each $\vec{n}$, if we allow repeats in the multi-indexes $\vec{n}$. The flexibility afforded by the presence of the scaling $d_{\vec{n}}$, and the allowance of repeats in the multi-indexes $\vec{n}$, leads immediately to a countably infinite sum of independent terms, but with new values of the $d_{\vec{n}}$'s, and possibly more than one monomial term per $\vec{n}$. The details are left to the reader.  

If such repeats are not allowed, the fact that $d_{\vec{n}}$ can be any analytic function of $|\vec{\theta}|^2$ does represent allowing correlations between polynomials of different order that differ only by powers of $|\vec{\theta}|^2$, but correlations between terms that have different angular behavior are forbidden if each $\vec{n}$ is assigned a single random coefficient. Some other treatments might instead expand in an orthonormal basis of functions, e.g. \cite{Mak24}. The basis functions are effectively the functions which are multiplied by the random coefficients to construct the emulator. Our functions here are simple products of various parameter components (monomials) multiplied by functions of $|\vec{\theta}|^2$, and are clearly not orthonormal. However, even with orthonormal terms one would be making some assumption that the terms are independent if each basis function is assigned a single random coefficient. It is not clear that orthonormality plays an important role in our application, but it is important to understand completeness. Our basis functions will become better defined below. The appendix, Sec. \ref{sec:basisprops}, reviews the orthonormality and completeness properties of the basis functions used here.

The dependence on $|\vec{\theta}|^2$ in the factors $d_{\vec{n}}$ and $f_K$ is not superfluous. As mentioned, these terms provide the flexibilty of correlating random coefficients for terms that differ only by factors of $|\vec{\theta}|^2$. However, the $|\vec{\theta}|$-dependence is necessary if the variance is to be independent of the choice in positioning the origin. At the origin, only the first term, $\vec{n}=0$, contributes to the emulator or to its variance. If there were no $|\vec{\theta}|$ dependence, the variance away from the origin would have the same contribution from the $\vec{n}=0$ term, but would also have contributions from other terms. Thus, to have a variance that is independent of position, the factors $d_{\vec{n}}$ and $f_K$ must include some behavior that decreases with $|\vec{\theta}|$ so that as $|\vec{\theta}|$ increases and terms of higher $K$ contribute more, terms of lower $K$ contribute less. In the next subsection it will be shown that this dependence must be a Gaussian form.

\subsection{Rotational and Translational Invariance}

The third assumption will be that the fluctuation at some point $\vec{\theta}$ is independent of the direction of $\vec{\theta}$. The fluctuation is
\begin{eqnarray}
\Var_0(Y(\vec{\theta})) &=&
\sum_{\vec{n}} 
f_K^2(|\vec{\theta}|^2)\sigma^2_A d_{\vec{n}}^2
\left(\frac{\theta_1^{2n_1}}{\Lambda^2}\right)\left(\frac{\theta_2^{2n_2}}{\Lambda^2}\right)\cdots 
\left(\frac{\theta_N^{2n_N}}{\Lambda^2}\right).
\end{eqnarray}
To be independent of direction, for a fixed norm the $\theta$ dependence must be proportional to $|\vec{\theta}|^{2K}$. Using the fact that
\begin{eqnarray}
|\vec{\theta}|^{2K}&=&\sum_{\vec{n},s.t. \sum_i n_i=K}\frac{K!}{n_1!\cdots n_N!}
(\theta_1^2)^{n_1}(\theta_2^2)^{n_2}\cdots (\theta_{N{\rm pars}}^2)^{n_{N{\rm pars}}},
\end{eqnarray}
one can recognize that
\begin{eqnarray}
d_{\vec{n}}^2(|\vec{\theta}|^2)\propto \frac{1}{n_1!n_2!\cdots n_{N_{\rm pars}}!}.
\end{eqnarray}
It can also be some arbitrary function of $K$ and $|\vec{\theta}|^2$, but that arbitrariness can be captured by $f_K(|\vec{\theta}|^2)$. Thus, with complete generality, once the independence of direction is enforced one can set
\begin{eqnarray}
d_{\vec{n}}(|\vec{\theta}|^2)&=& \frac{1}{\sqrt{n_1!n_2!\cdots n_{N_{\rm pars}}!}}.
\end{eqnarray}
The number of ways one can arrange the $K$ elements for a given $\vec{n}$ is $K!d_{\vec{n}}^2$.
The general form for the emulator, given the constraint that it is agnostic in regards to rotation in $\vec{\theta}$ space, is now
\begin{eqnarray}\label{eq:generalform}
Y(\vec{\theta})&=&\sum_{\vec{n}} 
f_{K(\vec{n})}(|\vec{\theta}|^2)
\frac{A_{\vec{n}}}{\sqrt{n_1!n_2!\cdots n_{N{\rm pars}}!}}
\left(\frac{\theta_1}{\Lambda}\right)^{n_1}\left(\frac{\theta_2}{\Lambda}\right)^{n_2}\cdots 
\left(\frac{\theta_N}{\Lambda}\right)^{n_{N{\rm pars}}}.
\end{eqnarray}

The fourth assumption is that the covariance for two different points should depend only on $|\vec{\theta}_1-\vec{\theta}_2|^2$. This will fix the form for $f_K(|\vec{\theta}|^2)$. The covariance is
\begin{eqnarray}\label{eq:EaEbgeneral}
\Cov_0(Y(\vec{\theta}_a),Y(\vec{\theta}_b))&=&\sigma_A^2
\sum_{\vec{n}}\frac{1}{\Lambda^{2K}}\frac{f_K(|\vec{\theta}_a|^2)f_K(|\vec{\theta}_b|^2)}{n_1!n_2!\cdots n_{N{\rm pars}}!}
(\theta_{a,1}\theta_{b,1})^{n_1}(\theta_{a,2}\theta_{b,2})^{n_2}\cdots(\theta_{a,N_{\rm pars}}\theta_{b,N_{\rm pars}})^{n_{N_{\rm pars}}}\\
\nonumber
&=&\sigma_A^2\sum_{\vec{n}}\frac{f_K(|\vec{\theta}_a|^2)f_K(|\vec{\theta}_b)|^2}{K!}(\vec{\theta}_a\cdot\vec{\theta}_b/\Lambda^2)^K.
\end{eqnarray}
The last step used the identity,
\begin{eqnarray}\label{eq:theta1dottheta2N}
|\vec{\theta}_a\cdot\vec{\theta}_b|^{K}&=&\sum_{\vec{n},s.t. \sum_i n_i=K}\frac{K!}{n_1!\cdots n_N!}
(\theta_{a,1}\theta_{b,1})^{n_1}(\theta_{a,2}\theta_{b,2})^{n_2}\cdots (\theta_{a,N{\rm pars}}\theta_{b,N{\rm pars}})^{n_{N{\rm pars}}},
\end{eqnarray}

If one were to choose
\begin{eqnarray}
f_K(|\vec{\theta}|^2)&=e^{-|\vec{\theta}|^2/2\Lambda^2},
\end{eqnarray}
the resulting covariance would be
\begin{eqnarray}
\Cov_0(Y(\vec{\theta}_a),Y(\vec{\theta}_b))&=&\sigma_A^2e^{-|\vec{\theta}_a-\vec{\theta}_b|^2/2\Lambda^2},
\end{eqnarray}
which would satisfy the constraint of being a function of $|\vec{\theta}_a-\vec{\theta}_b|$. But, one can go further and show that this is the only such solution. To demonstrate this, one can write 
\begin{eqnarray}\label{eq:fKdefault}
f_K(|\vec{\theta}|^2)&=&e^{-\frac{|\vec{\theta}|^2}{2\Lambda^2}}+\delta f_K(|\vec{\theta}|^2).
\end{eqnarray}
Further, one can expand $\delta f_K$,
\begin{eqnarray}
\delta f_K(|\vec{\theta}|^2)=\alpha_{K,2}|\vec{\theta}|^4+\alpha_{K,3}|\vec{\theta}|^6+\alpha_{K,4}|\vec{\theta}|^8+\cdots
\end{eqnarray}
The series begins for quartic terms because adding a quadratic term, beyond what is in the Gaussian in Eq. (\ref{eq:fKdefault}), would violate the constraint that overall the quadratic term in Eq. (\ref{eq:EaEbgeneral}) must behave as $|\vec{\theta}_a-\vec{\theta}_b|^2$. The corrections due to $\delta f_K$ of quartic behavior are
\begin{eqnarray}
\delta \Cov_0( Y(\vec{\theta}_a),Y(\vec{\theta}_b))&=&
\alpha_{0,2}|\vec{\theta}_a|^4+\alpha_{0,2}|\vec{\theta}_b|^4+\alpha_{2,0}(\vec{\theta}_a\cdot\vec{\theta}_b)^2.
\end{eqnarray}
One can set $\alpha_{0,2}=-\alpha_{2,0}/2$ to make this disappear when $\vec{\theta}_1=\vec{\theta}_2$, but this correction is not proportional to $|\vec{\theta}_1-\vec{\theta}_2|^4$. The same can be said for higher-order terms in $\theta$. Thus, $\delta f$ must be zero and the form for $f_K(|\vec{\theta}|^2)$ is completely determined to match that of Eq. (\ref{eq:fKdefault}) if one insists that the covariance be a function of $|\vec{\theta}_1-\vec{\theta}_2|^2$.

After enforcing the four criteria, the form for the emulator and covariance are now fixed to
\begin{eqnarray}\label{eq:E2default}
Y(\vec{\theta})&=&e^{-|\vec{\theta}|^2/2\Lambda^2}\sum_{\vec{n}} 
\frac{A_{\vec{n}}}{\sqrt{n_1!n_2!\cdots n_{N{\rm pars}}!}}
\left(\frac{\theta_1}{\Lambda}\right)^{n_1}\left(\frac{\theta_2}{\Lambda}\right)^{n_2}\cdots 
\left(\frac{\theta_N}{\Lambda}\right)^{n_{N{\rm pars}}},\\
\nonumber
\Cov_0(Y(\vec{\theta}_a),Y(\vec{\theta}_b))&=&\sigma_A^2e^{-|\vec{\theta}_a-\vec{\theta}_b|^2/2\Lambda^2},\\
\nonumber
\Var_0(Y(\vec{\theta}))&=&\sigma_A^2.
\end{eqnarray}

In fact, for the four criteria mentioned above, satisfying the third criteria is automatically satisfied if the fourth is satisfied, because the fluctuation is simply the covariance evaluated at $\vec{\theta}_1=\vec{\theta}_2$. Thus, only three criteria are required to constrain the form to that of Eq. (\ref{eq:E2default}). The form for the emulator and covariances is fixed to that of Eq. (\ref{eq:E2default}) if the following choices or assumptions are made:
\begin{enumerate}\itemsep=0pt
\item The function should be analytic throughout the prior.
\item At any point the function should vary according to the same Gaussian distribution.
\item The prior-to-constraint covariance, $\Cov_0(Y(\vec{\theta}_a),Y(\vec{\theta}_b))$, should depend only on $|\vec{\theta}_1-\vec{\theta}_2|^2$.
\end{enumerate}
From here on, the covariance will be referred through the function $C_0(\vec{\theta}_a,\vec{\theta}_b)$, defined as
\begin{eqnarray}\label{eq:C0default}
C_0(\vec{\theta}_a,\vec{\theta}_b)&\equiv&\frac{1}{\sigma_A^2}\Cov_0(Y(\vec{\theta}_a),Y(\vec{\theta}_b))=e^{-|\vec{\theta}_a-\vec{\theta}_b|^2/2\Lambda^2}.
\end{eqnarray}
The subscript zero again refers to the fact that this averaging is performed before one has applied any constraints from the training data.

\subsection{Higher Moments and the Assumption that Coefficients have Gaussian Distributions}
\label{sec:Agauss}

One might ask whether all other moments, e.g. $\E_0(Y^4(\vec{\theta}))$, are also independent of $\vec{\theta}$. The answer is yes, and derives from the initial choice that the weights for the coefficients are Gaussian. For Gaussians centered at zero, all cumulants are zero except for the second cumulant, which is $\sigma_A^2$. Thus, when averaging over the initial coefficients, the averaged fourth moment of the coefficients becomes 
\begin{eqnarray}\label{eq:Acumulant}
\Cov_0(A_i,A_j,A_m,A_n)&=&
\Cov_0(A_i,A_j)\Cov_0(A_m,A_n)\\
\nonumber&&
+\Cov_0(A_i,A_m)\Cov_0(A_j,A_n)
+\Cov_0(A_i,A_n)\Cov_0(A_j,A_m).
\end{eqnarray}
For any distribution where $\E_0(A_i)=0$, the independent coefficients, Eq. (\ref{eq:Acumulant}) would be zero if the indices were not all the same, e.g. if $i=j$ and $m=n$, but $i\ne m$. Gaussian distributions are special in that this holds true even if $i=j=m=n$. If one expands $Y(\vec{\theta})$ in terms of the coefficients and polynomials, one readily sees that this implies that
\begin{eqnarray}
\Cov_0(Y(\vec{\theta}_a),Y(\vec{\theta}_b),Y(\vec{\theta}_c),Y(\vec{\theta}_c))&=&
\Cov_0(Y(\vec{\theta}_a),Y(\vec{\theta}_b))\Cov_0(Y(\vec{\theta}_c),Y(\vec{\theta}_d))\\
\nonumber
&+&\Cov_0(Y(\vec{\theta}_a),Y(\vec{\theta}_c))\Cov_0(Y(\vec{\theta}_b),Y(\vec{\theta}_d))
+\Cov_0(Y(\vec{\theta}_a),Y(\vec{\theta}_d))\Cov_0(Y(\vec{\theta}_b),Y(\vec{\theta}_c)).
\end{eqnarray}
All covariances to order $Y^4$ can therefore be expressed in terms of covariances of order $Y^2$. If the covariances of order $Y^2$ are independent of translations and rotations, the covariances of order $Y^4$ must then also satisfy this requirement.  If all the coordinates of a covariance are rotated or translated, the covariance is unchanged. For the particular case where one considers all four angles to be equal to one another this gives
\begin{eqnarray}
\E_0(Y^4(\vec{\theta}))&=&3\sigma_A^4,
\end{eqnarray}
which is independent of $\vec{\theta}$. In fact, because all cumulants disappear, covariances to any order, $\Cov_0(Y(\vec{\theta}_1,\cdots Y(\vec{\theta}_n))$, can be expressed in terms of $\Cov_0(Y(\vec{\theta}_i),Y(\vec{\theta}_j))=\sigma_A^2e^{-|\vec{\theta}_i-\vec{\theta}_j|^2/2\Lambda^2}$.

If one replaces the choice of Gaussian distributions for the coefficients with another distribution with the same second moments, with all the odd moments being zero, but with different cumulants for the higher even moments, the property that the covariances will only depend on relative coordinates will be broken. For example, if one considers the fourth moment,
\begin{eqnarray}
    \E_0(A_i^4)&=&3\sigma_A^4+\kappa.
\end{eqnarray}
Here, $\kappa=0$ for a Gaussian distribution. If one were to calculate the fourth moment of $E$ with all four positions being equal, and with $\vec{\theta}$ being along the first axis, $\vec{\theta}=[\theta_1,0,0\cdots]$,
\begin{eqnarray}
    \E_0(Y(\vec{\theta})^4)&=&3\sigma_A^2+\kappa e^{-2\theta_1^2|^2/\Lambda^2}
    \sum_{n}\frac{1}{(n!)^2}\theta_1^{4n}.
\end{eqnarray}
Clearly, this is not independent of $\theta_1$, and is clearly not a function of relative coordinates, as all the relative coordinates are zero. The correlations are then clearly not invariant under translation and rotation. Thus, if one wishes all the higher-order correlations to depend only on relative coordinates and be independent of translation and rotation, the distribution for the coefficients must have zero cumulants. The distribution $W(A_{\vec{n}})$ in Eq. (\ref{eq:AGauss}) should therefore be of a Gaussian form.

\subsection{Comparison with Gaussian Process Emulators}

Ultimately, an emulator constructed from the form proposed above depends only on the unconstrained correlations, $C_0(\vec{\theta}_i,\vec{\theta}_j)$, a result known as the kernel trick. The emulators then map perfectly to Gaussian process emulators, which are also perfectly defined once one knows $C_0(\vec{\theta}_i,\vec{\theta}_j)$. Here, the convergence parameter $\Lambda$ is introduced to set the relative strength of order $K+1$ terms to order $K$ terms, whereas in Gaussian process emulators the same parameter is introduced to set the correlation of the unconstrained correlation. Despite these differences in starting points and perspective, the resulting emulators are identical. This equivalence demonstrates how Gaussian process emulators can be equivalently viewed as polynomial expansions if one uses the basis assumption in Eq. (\ref{eq:generalform0}).

Gaussian process emulators often assume the correlation has a more general form,
\begin{eqnarray}\label{eq:GPform}
C_0(\vec{\theta}_1,\vec{\theta}_2)&\sim&\exp\left\{\left(-\frac{|\vec{\theta}_1-\vec{\theta}_2|^2}{2\Lambda^2}\right)^{\alpha},\right\},
\end{eqnarray}
where $\alpha$ is an additionally adjustable hyper parameter. For $\alpha=1$ this corresponds precisely to the form discussed above. If the Taylor expansion for the emulator is analytic, one should  expect that the kernel, $C_0(\vec{\theta}_1,\vec{\theta}_2)$, should also be analytic in terms of $\vec{\theta}$. The form expressed in Eq. (\ref{eq:GPform}) will only be analytic if $\alpha$ is an integer. Thus, the Gaussian form for the emulator in Eq. (\ref{eq:E2default}) should not come as a surprise. What is perhaps surprising is that $\alpha$ cannot be some other integer. 

Despite the fact that the emulators are identical in practice to GP forms, there is insight into viewing the emulator from the perspective of polynomial expansions. Once one considers the emulator form in terms of terms of a polynomial expansion, one also understands how many training points are needed to constrain all terms of a specific order. For example, for a linear fit, $N_{\rm pars}+1$ training points are needed to determine the intercept and slopes. For a quadratic fit, one would need $(N_{\rm pars}+1)(N_{\rm pars}+2)/2$ terms, and for a cubic fit one would need $(N_{\rm pars}+1)(N_{\rm pars}+2)(N_{\rm pars}+3)/6$ training points. Thus, it might motivate one to choose a number of training points to constrain all coefficients to a given order. 

\subsection{Tuning the Emulator}

Due to the Gaussian weights, the most probable set of expansion coefficients would be when all $A_{\vec{n}}$ are zero. However, once one fits to $N_{\rm train}$ training points one needs to find the set of coefficients that both maximizes the weight function, $W(A_{\vec{n}})$, while reproducing the $N_{\rm train}$ training points. Whereas there might be dozens or a few hundred training points, there are an infinite number of coefficients. Thus, one needs to maximize the weight with $N_{\rm train}$ constraints. Here, it is shown how one can analytically solve for the set of most probable coefficients, $A_{\vec{n}}$, given the training values, $y_a=y(\vec{\theta}_a)$, where $\vec{\theta}_a$ are the points at which the emulator was trained. For shorthand, we define the function $T_{\vec{n}}(\vec{\theta})$ to rewrite the emulator in Eq. (\ref{eq:E2default}),
\begin{eqnarray}\label{eq:ETdef}
Y(\vec{\theta})&=&\sum_{\vec{n}}A_{\vec{n}}\mathcal{T}_{\vec{n}}(\vec{\theta}),\\
\nonumber
\mathcal{T}_{\vec{n}}(\vec{\theta})&=&
\frac{1}{\sqrt{n_1!n_2!\cdots n_{N{\rm pars}}!}}
\left(\frac{\theta_1}{\Lambda}\right)^{n_1}\left(\frac{\theta_2}{\Lambda}\right)^{n_2}\cdots 
\left(\frac{\theta_N}{\Lambda}\right)^{n_{N{\rm pars}}}e^{-|\vec{\theta}|^2/2\Lambda^2}.
\end{eqnarray}
Here, $T_{\vec{n}}(\vec{\theta})$ is shorthand for the basis functions used to represent the emulator. None of the following derivations are dependent on this particular set of basis functions.

Evaluated at the $N_{\rm train}$ training points, one can define the $N_{\rm train}\times N_{\rm coef}$ matrix,
\begin{eqnarray}
T_{a,\vec{n}}\equiv\mathcal{T}_{\vec{n}}(\vec{\theta}_a),
\end{eqnarray}
where $\vec{\theta}_a$ labels the training points, $0< a< N_{\rm train}$. The number of coefficients, $N_{\rm coef}$, approaches infinity.

The goal is to find the coefficients, $A_{\vec{n}}$, that maximize the probability given the constraints of matching the training data. If the training data are to be exactly matched, the differential probability of a given set of coefficients is
\begin{eqnarray}\label{eq:Nexact}
dP(\vec{A},\sigma_A,\Lambda)&\sim &
\prod_{\vec{n}}\frac{dA_{\vec{n}}}{(2\pi\sigma_A^2)^{1/2}}
e^{-A_{\vec{n}}^2/2\sigma_A^2}
\prod_{a=1}^{ N_{\rm train}}\delta[y(\vec{\theta}_a)-Y(\vec{\theta}_a)].
\end{eqnarray}
If the training constraints are not exact, the delta function should be replaced with some other form, perhaps a Gaussian. The Gaussian would be reasonable if the training had some sort of point-by-point noise, such as one would have if the full-model calculation involved some sort of sampling or Monte Carlo evaluations. However, more generally, one could also consider the fact that the training calculation might have some sort of error, e.g. that of missing physics. In that cast the error might be correlated between points.  Noise in the full-model calculations at points $\vec{\theta}_a$ and $\vec{\theta}_b$ can be characterized by some error matrix, 
\begin{eqnarray}
\Cov(y_a,y_b)&=&\sigma_A^2\Delta_{ab}.
\end{eqnarray}
Here, $y_a=y(\vec{\theta}_a)$, the full model value calculated at the training point $a$, and $\delta y_a$ is the uncertainty of the full model due to point-by-point noise in the full model or any other source of error.  The choice to scale out the factor of $\sigma_A^2$ in the definition is so that $\Delta_{ab}$ is dimensionless. If the random noise is completely uncorrelated, $\Delta_{ab}$ is diagonal. For example, if one were to recalculate the full model at exactly the same value of $\vec{\theta}_a$, the new value of $y_a$ could be different. This is often the case for models involving Monte Carlo calculation or if numerical noise is non-negligible. Given that the characteristic size of the observable is $\sigma_A$, the elements of the diagonal matrix would then be chosen according to the percentage random error. For example, if the random error were 1\%, one would chose $\Delta_{ab}=\delta_{ab}10^{-4}$.  

To account for the uncertainty represented by $\Delta$, one can relax the constraint imposed by the delta function by replacing it, 
\begin{eqnarray}\label{eq:deltareplace1}
\prod_a\delta[y_a-Y(\vec{\theta}_a)]&\rightarrow &
\frac{|\Delta|^{-1/2}}{(2\pi)^{N_{\rm train}/2}}\exp\left\{\frac{-1}{2}
(y_a-Y(\vec{\theta}_a))\Delta_{ab}^{-1}(y_b-Y(\vec{\theta}_b))
\right\}.
\end{eqnarray}
The second step above exploits the fact that the Fourier transform of a Gaussian is a Gaussian. The r.h.s. above behaves like a delta function in the limit that the eigenvalues of $\Delta$ approach infinity. Otherwise, in the basis where $\Delta$ is diagonal this can be thought of as a finite-width modification to a delta function. Next, one can replace $Y(\vec{\theta})$ in Eq. (\ref{eq:deltareplace1}) using Eq. (\ref{eq:ETdef}). The altered $\delta$ functions then become
\begin{eqnarray}
\prod_a\delta[y_a-Y(\vec{\theta}_a)]&\rightarrow &\int\prod_a\left(\frac{dk_a}{2\pi}\right)\exp\left\{ik_a(y_a-\sum_{\vec{m}}A_{\vec{m}}T_{\vec{m}}(\vec{\theta}_a))-\sum_{b}\frac{\sigma_A^2}{2}k_a\Delta_{ab}k_b\right\}.
\end{eqnarray}
Defining the net probability after integrating over the coefficients above, but not over $\sigma_A^2$ or $\Lambda$,
\begin{eqnarray}
\label{eq:Zdef}
\mathcal{Z}(\sigma_A^2,\Lambda)&\equiv&\int\left(\prod_{\vec{n}}dA_{\vec{n}}\right)~\frac{dP(\vec{A},\sigma_A,\Lambda)}{\prod dA_{\vec{n}}}\\
\nonumber
&=&\int \left(\prod_{\vec{n}} \frac{dA_{\vec{n}}}{\sqrt{2\pi\sigma_A^2}}\right)
~\exp\left[
-\frac{\sum_iA_{\vec{i}}^2}{2\sigma_A^2}\right]\int \frac{dk}{2\pi}\exp\left[i\sum_bk_b(y_b-\sum_\ell A_\ell T_\ell(\vec{\theta}_b))-\frac{\sigma_A^2}{2}\sum_{ab}k_a\Delta_{ab}k_b
\right].
\end{eqnarray}
The function $Z(\sigma_A^2,\Lambda)$ is the conditional probability of observing the training data given $\sigma_A$ and $\Lambda$ given the training data. From Eq. (\ref{eq:Nexact}) one can see that $Z(\sigma_A^2,\Lambda)$ is normalized when integrated over all possible training data,
\begin{equation}\label{eq:Znorm}
\int\left(\prod_a dy_a\right) Z(\sigma_A^2,\Lambda)=1.
\end{equation}
The average coefficients, and correlation of coefficients, for a given $\sigma_A$ and $\Lambda$, become 
\begin{eqnarray}
\E(A_{\vec{n}})&=& \mathcal{Z}^{-1}
\int \left(
\prod_{\vec{i}}
\frac{dA_{\vec{i}}}{\sqrt{2\pi\sigma_A^2}}
\right)~
\exp\left[
-\sum_{\vec{j}} \frac{A_{\vec{j}}^2}{2\sigma_A^2}
\right]
A_{\vec{n}}\\
\nonumber 
&&
\int \left(
\prod_a\frac{dk_a}{2\pi}
\right)
\exp\left[i\sum_bk_b(y_b-\sum_{\vec{\ell}}A_{\vec{\ell}}T_{\vec{\ell}}(\vec{\theta}_b))-\frac{\sigma_A^2}{2}\sum_{ab}k_a\Delta_{ab}k_b\right],\\
\nonumber
\Cov(A_{\vec{m}},A_{\vec{n}})&=& \mathcal{Z}^{-1}\int \left(\prod_{\vec{i}} \frac{dA_{\vec{i}}}{\sqrt{2\pi\sigma_A^2}}\right)~  \exp\left[
-\sum_j\frac{A_{\vec{j}}^2}{2\sigma_A^2}-\frac{\sigma_A^2}{2}\sum_{ab}k_a\Delta_{ab}k_b\right]\\
\nonumber
&&
A_{\vec{m}}A_{\vec{n}}\int \prod_a\frac{dk_a}{2\pi}\exp\left[i\sum_bk_b(y_b-\sum_{\vec{\ell}}A_{\vec{\ell}}T_{\vec{\ell}}(\vec{\theta}_b))\right],
\end{eqnarray}
In this subsection, the means, variances and covariances refer to probabilities for assumed values of the hyper-parameters, $\sigma_A$ and $\Lambda$.

After completing the square, the covariance can be written as
\begin{eqnarray}
\Cov(A_{\vec{m}},A_{\vec{n}})&=& \mathcal{Z}^{-1}\int\left(\prod_a\frac{dk_a}{2\pi}\right)
e^{i\sum_bk_by_b-\frac{\sigma_A^2}{2}\sum_{ab}k_a\Delta_{ab}k_b}~\\
\nonumber
&&~~\int \left(\prod_{\vec{i}} dA_{\vec{i}}\right)~A_{\vec{m}}A_{\vec{n}}
\exp\left[-\frac{1}{2\sigma_A^2}\left(A_{\vec{i}}+i\sigma_A^2\sum_bT_{\vec{i}}](\vec{\theta}_b)k_b\right)^2
-\frac{\sigma_A^2}{2}\left(\sum_{b,\vec{\ell}} T_{\vec{\ell}}(\vec{\theta}_b)k_b\right)^2
\right]\\
\nonumber
&=& \mathcal{Z}^{-1}
\left(\int\prod_a\frac{dk_a}{2\pi}\right)~\exp\left[-\frac{\sigma_A^2}{2}(T_{\vec{m}}(\vec{\theta}_b)k_b)^2+ik_by_b-\frac{\sigma_A^2}{2}k_a\Delta_{ab}k_b\right]\\
\nonumber
&&~~\left(\int \prod_{\vec{i}} dA_{\vec{i}}\right)~(A_{\vec{n}}-i\sigma_A^2\sum_bT_k(\vec{\theta}_b)k_b)(A_{\vec{m}}-i\sigma_A^2T_{\vec{m}}(\vec{\theta}_c)k_c)
\exp\left[-\sum_{\vec{\ell}}\frac{A_{\vec{\ell}}^2}{2\sigma_A^2}\right]
\end{eqnarray}
Similarly,
\begin{eqnarray}
E(A_{\vec{m}}) &=&-\mathcal{Z}^{-1}
\int\left(\prod_a\frac{dk_a}{2\pi}\right)~\exp\left[-\sum_b\frac{\sigma_A^2}{2}(T_{\vec{m}}(\vec{\theta}_b)k_b)^2+i\sum_bk_by_b-\frac{\sigma_A^2}{2}\sum_{ab}k_a\Delta_{ab}k_b\right]\\
\nonumber
&&~~\int \left(\prod_{\vec{i}} \frac{dA_{\vec{i}}}{\sqrt{2\pi\sigma_A^2}}\right)~i\sigma_A^2\left(\sum_cT_{\vec{m}}(\vec{\theta}_c)k_c\right)
\exp\left[-\sum_{\vec{\ell}}\frac{A_{\vec{\ell}}^2}{2\sigma_A^2}\right].
\end{eqnarray}

It is convenient to define the function,
\begin{eqnarray}
C_0(\vec{\theta}_1,\vec{\theta}_2)&\equiv&\sum_{\vec{i}} T_{\vec{i}}(\vec{\theta}_1)T_{\vec{i}}(\vec{\theta}_2)\\
\nonumber
&=&\frac{\Cov_0(Y(\vec{\theta}_1),Y(\vec{\theta}_2))}{\sigma_A^2}.
\end{eqnarray}
The subscript zero, denotes that this is the average in the absence of training. 

It is convenient to define the $N_{\rm train}\times N_{\rm train}$ matrix,
\begin{eqnarray}\label{eq:Bdef}
 B_{ab}&\equiv& C_0(\vec{\theta}_a,\vec{\theta}_b)+\Delta_{ab},
\end{eqnarray}
where $\vec{\theta}_a$ are the locations of the training points. Performing the integrals for $ \mathcal{Z}(\sigma_A^2,\Lambda)$ in Eq. (\ref{eq:Zdef}),
\begin{eqnarray}\label{eq:netexactprob}
 \mathcal{Z}(\sigma_A^2,\Lambda)&=&\sigma_A^{-N_{\rm train}}\left({\rm det} B\right)^{-1/2}\exp\left[-\frac{1}{2\sigma_A^2}\sum_{ab}y_a B^{-1}_{ab}y_b\right].
\end{eqnarray}
One can then find the moments of the integral over $k$,
\begin{eqnarray}
\E(k_a)&=&-i\partial_{y_a}\ln( \mathcal{Z})=\frac{i}{\sigma_A^2} \sum_bB^{-1}_{ab}y_b,\\
\nonumber
\Cov(k_a,k_b)&=&-\partial_{y_a}\partial_{y_b}\ln( \mathcal{Z})=\frac{1}{\sigma_A^2} B_{ab}^{-1}.
\end{eqnarray}
The averaged coefficients are then:
\begin{eqnarray}\label{eq:exactmoments}
\E(A_{\vec{n}})&=&\sum_{ab}T_{\vec{n}}(\vec{\theta}_a) B^{-1}_{ab}y_b,\\
\nonumber
\Cov(A_{\vec{m}},A_{\vec{n}})&=&
\delta_{\vec{m},\vec{n}}\sigma_A^2-\sigma_A^2\sum_{ab} 
T_{\vec{m}}(\vec{\theta}_a)B^{-1}_{ab}T_{\vec{n}}(\vec{\theta_b}).
\end{eqnarray}

The observables and their fluctuations due to the fluctuation in the coefficient are then:
\begin{eqnarray}
\E(Y(\vec{\theta}))&=&\sum_{\vec{n}}\E(A_{\vec{n}})T_{\vec{n}}(\vec{\theta}),\\
\nonumber
&=&\sum_{ab}y_aB^{-1}_{ab}C_0(\vec{\theta}_b,\vec{\theta}_a),\\
\nonumber
\Var(Y(\vec{\theta}))&=&
\sum_{\vec{m}\vec{n}}\Cov(A_{\vec{m}},A_{\vec{n}}) T_{\vec{m}}(\vec{\theta})T_{\vec{n}}(\vec{\theta})\\
\nonumber
&=&\sigma_A^2\left[
C_0(\vec{\theta},\vec{\theta})-\sum_{ab}C_0(\vec{\theta},\vec{\theta}_a)B^{-1}_{ab}C_0(\vec{\theta}_b,\vec{\theta})
\right].
\end{eqnarray}
The fact that the prediction for  $Y(\vec{\theta})$ and its uncertainty can be expressed so simply in terms of the correlation kernel $C_0$ without involving the coefficients is the well known kernel trick \cite{RasmussenGP}. This enables one to build the emulator without every bothering to calculate the individual expansion coefficients $A_{\vec{n}}$

After defining
\begin{eqnarray}
\chi_a&\equiv \sum_bB^{-1}_{ab}y_b,
\end{eqnarray}
the emulated value and emulator uncertainty at some point $\vec{\theta}$ are
\begin{eqnarray}
\label{eq:exactydely2}
Y(\vec{\theta})&=&\sum_a\chi_aC_0(\vec{\theta}_a,\vec{\theta}),\\
\nonumber
\Var(Y(\vec{\theta}))&=&\sigma^2_Y(\vec{\theta})=
\sigma_A^2 C_0(\vec{\theta},\vec{\theta})
-\sigma_A^2\sum_{ab}C_0(\vec{\theta},\vec{\theta}_a) B^{-1}_{ab}C_0(\vec{\theta}_b,\vec{\theta}).
\end{eqnarray}
The optimum value $\E(Y(\vec{\theta}))$ depends linearly on the training values, $y_{a=1,N_{\rm train}}$. The emulator's uncertainty depends on $\vec{\theta}$ and vanishes at the training points if $\Delta_{ab}=0$. It is linear in $\sigma_A$ and has a more complicated dependence on $\Lambda$. It does not depend on the model values at those training points, $y_a$, although after tuning the choice of $\sigma_A$ will depend on the training values. Not surprisingly, this is the exactly the same expression used for Gaussian process emulators, with the requirement that the form for the GP correlation kernel should be Gaussian.

Even though the procedure outlined above allows one to solve for all the coefficients, $A_{\vec{n}}$, solving for the individual coefficients is not necessary because both the emulator and its uncertainty can be expressed, as shown in Eq. (\ref{eq:exactydely2}), in terms of the bare correlation kernel, $C_0(\vec{\theta}_1,\vec{\theta}_2)$. The simple form is the same for some GP emulators, but rather than being arbitrarily assigned as is the case for the GP paradigm, in this perspective the form is motivated by the three simple assumptions: analyticity, independent Gaussian fluctuation of the polynomial terms, and invariance to translations and rotations. The fact that the emulator can be expressed simply in terms of the kernel, $C_0(\vec{\theta}_1,\vec{\theta}_2)$, is not unique to the polynomial expansion used here. Any form where $Y(\vec{\theta})$ is expressed as a sum of terms, where each term is linear in some coefficient, $A_i$, that has some Gaussian distribution, can be expressed in the form of Eq. (\ref{eq:exactydely2}). This is known as the kernel trick \cite{RasmussenGP} and has been exploited in many previous studies involving polynomial expansions. For example, if one had chosen the form $f_K(|\vec{\theta}|^2)=1$, rather than being a Gaussian, the fluctuation would then not be independent of the norm of $\vec{\theta}$, and the correlation would no longer depend only on $|\vec{\theta}_1-\vec{\theta}_2|^2$, but one could still optimize the expansion. The resulting correlation kernel would change to $C_0(\vec{\theta}_1,\vec{\theta}_2)=e^{\vec{\theta}_1\cdot\vec{\theta}_2/\Lambda^2}$, but the form of Eq. (\ref{eq:exactydely2}) would not change.

Despite the fact that the emulator's behavior is completely determined by the kernel, knowing the form for the polynomial expansion that satisfies the three aforementioned constraints is nonetheless useful for testing. In Sec. \ref{sec:testing} random test functions will be created using this form. Each function is characterized by $\sigma_A,\Lambda$ and the set of random Gaussian coefficients. If one trains an emulator on such a test function, the emulator's value of $\sigma_A$ and $\Lambda$ should be close to those used to generate the function. Such tests are presented in Sec. \ref{sec:testing}.

The fluctuation described in Eq. (\ref{eq:exactydely2}) were only those due to the fluctuation of the coefficients. This fluctuation was calculated ignoring the fact that the hyper-parameters $\sigma_A$ and $\Lambda$ are not precisely determined. An additional contribution to $\Var(Y(\vec{\theta}))$ arises if one accounts for the variance in $\Lambda$. This is the subject of the next subsection.

\subsection{Finding the Optimized Hyper-Parameters}\label{sec:choosehyper}

This approach to emulation has two hyper parameters, the variance of the coefficients, $\sigma_A$, and the convergence parameter $\Lambda$. If one has sufficient training data, one can find their optimum values. The next subsection considers choosing $\sigma_A$ for some fixed $\Lambda$ and the following subsection considers choosing $\Lambda$ while fixing $\sigma_A$. 

Both $\sigma_A$ and $\Lambda$ are unknown, and the probability distribution for them is determined by the function $Z(\sigma_A^2,\Lambda)$, which can be considered as the conditional probability of observing the training data given $\sigma_A$ and $\Lambda$. Restating Eq. (\ref{eq:netexactprob}),
\begin{eqnarray}\label{eq:ZofSigmaALambda}
\mathcal{Z}(\sigma_A^2,\Lambda)&=&\int \prod_{\vec{i}} dA_{\vec{i}} P(\vec{A},\sigma_A^2,\Lambda)\\
\nonumber
&=&\sigma_A^{-N_{\rm train}}\left({\rm det} B\right)^{-1/2}\exp\left[-\frac{1}{2\sigma_A^2}y_a B^{-1}_{ab}y_b\right].
\end{eqnarray} 
By Bayes theorem the probability of $\sigma_A$ and $\Lambda$ is proportional to the product of $Z$ and $Q(\sigma_A^2,\Lambda)$, the prior probability of $\sigma_A$ and $\Lambda$. I.e. $Q(\sigma_A^2,\Lambda)$ is the probability for $\sigma_A$ and $\Lambda$ in the absence of training data.
The most appealing method to chose $\sigma_A$ and $\Lambda$, would be to choose their average over the weight $ZQ$,
\begin{eqnarray}
    \bar{\sigma}_A &=&\frac{1}{\mathcal{N}}\int d\sigma_Ad\Lambda~\sigma_AZ(\sigma_A^2,\Lambda)Q(\sigma_A^2,\Lambda),\\
    \nonumber
    \bar{\Lambda} &=&\frac{1}{\mathcal{N}}\int d\sigma_Ad\Lambda~\Lambda ~Z(\sigma_A^2,\Lambda)Q(\sigma_A^2,\Lambda),\\
    \nonumber
    \mathcal{N}&=&\int d\sigma_A d\Lambda Z(\sigma_A^2,\Lambda)Q(\sigma_A^2,\Lambda).
\end{eqnarray}
Similarly, one can find expressions for the uncertainties of the hyper parameters by calculating higher-moment integrals. We define the covariance matrix  $W_{\alpha,\beta}$,
\begin{eqnarray}
    W_{\sigma_A\sigma_A}&=&\Var(\sigma_A),\\
    \nonumber
    W_{\sigma_A\Lambda}&=&W_{\Lambda\sigma_A}=\Cov(\sigma_A,\Lambda),\\
    \nonumber
    W_{\Lambda\Lambda}&=&\Var(\Lambda).
\end{eqnarray}
Here, the averages, variances and covariances take into account the weights $Z(\sigma_a^2,\Lambda)Q(\sigma_A^2,\Lambda)$.  Unfortunately, the functions $Z(\sigma_A,\Lambda)$ approaches a constant as $\Lambda\rightarrow\infty$. This means that all moments, even the normalization, are strongly affected by the range of the prior, $Q(\sigma_A^2,\Lambda)$. In fact, if the prior has some range $\Lambda_{\rm max}$, the $n^{\rm th}$ moment would scale as $\Lambda_{\rm max}^n$ for large $\Lambda_{\rm max}$. Thus, this weight gives one pause in using moments of the distribution to represent best values or uncertainties.

As an alternative to integrating the expressions above over $\Lambda$ in order to find the optimized values of $\sigma_A$ and $\Lambda$,  one can solve for the values of $\sigma_A$ and $\Lambda$ that maximize $Z(\sigma_A^2,\Lambda)$ as alternatives to $\bar{\Lambda}$ and $\bar{\sigma}_A$. Further, one can assign the uncertainty by calculating the curvatures, i.e. the second derivatives of $Z(\sigma_A^2,\Lambda)$ evaluated at the values of $\sigma_A$ and $\Lambda$ for which the first derivatives vanish. 
\begin{eqnarray}\label{eq:Walphabetaalt}
    W_{\alpha\beta}&=&
    -\frac{\partial^2 \ln(Z)}{\partial\alpha\partial\beta}.
\end{eqnarray}
Essentially, this corresponds to assuming that $\ln(Z)$ is quadratic in the parameters. Normally, choosing the point of the maximum and its curvature is a reasonable assumption whenever $Z$ is a sharply peaked function. If $Z$ were a Gaussian function of the parameters, which it is not, one could assume an arbitrarily wide prior, and the averages of parameters and their covariances would be equivalent to the position of the maximum and the curvature at the maximum.

Another advantage of using the alternative expressions for $\sigma_A$ and $\Lambda$ and their uncertainties are that the quantities are calculated  quickly. A Newtons method algorithm required $\lesssim 10$ evaluations of $Z(\sigma_A^2,\Lambda)$ and its derivatives. Both $Z$ and its derivatives are easily calculated and finding these points at which $Z$ is maximized is robust. The calculations shown in this paper apply this alternative. As a test, in Sec. \ref{sec:testing} functions $y(\vec{\theta})$ are built with polynomial coefficients chosen consistent with the form assumed by the emulator with specific values of $\sigma_A$ and $\Lambda$. Emulators, $Y(\vec{\theta})$ are then constructed with $\sigma_A$ and $\Lambda$ chosen according to the method described above. The chosen values of $\sigma_A$ and $\Lambda$ well match those used to construct $y(\vec{\theta})$ for larger numbers training points.

\subsubsection{Choosing the Characteristic Coefficient Width}
First, we consider fixing the convergence parameter $\Lambda$ and finding the optimum value of $\sigma_A$. This involves maximizing the probability $Z(\sigma_A^2,\Lambda)$, defined in Eq. (\ref{eq:Zdef}) for a fixed $\Lambda$. 
The goal here is to find the value of $\sigma_A$ that maximizes $ \mathcal{Z}$. For the case where $\Delta=0$ or if $\Delta$ is independent of $\sigma_A$, one can find $\sigma_A$ rather easily because $B_{ab}$ does not depend on $\sigma_A$. In that case, setting the derivative w.r.t. $\sigma_A^2$ equal to zero,
\begin{eqnarray}
0&=&\frac{N_{\rm train}}{\sigma_A^2} \mathcal{Z}-\frac{\mathcal{Z}}{\sigma_A^4}y_a B^{-1}_{ab}y_b,\\
\label{eq:sigmaAguess}
\sigma_A^2&=&\frac{1}{N_{\rm train}}y_a B^{-1}_{ab}y_b.
\end{eqnarray}

\subsubsection{Choosing the Convergence Parameter}

Eq. (\ref{eq:sigmaAguess}) provides the means for choosing $\sigma_A$ for a given value $\Lambda$. One can then plot $Z(\sigma_A^2,\Lambda)$, as given in Eq. (\ref{eq:ZofSigmaALambda}), as a function of $\Lambda$ and choose the value of $\Lambda$ that maximizes $Z$. For some value of $\Lambda$ the expressions for $B_{ab}$ and $\sigma_A$ are easily found. One can also find analytic expressions for the derivatives of $B$ with respect to $\Lambda$. This allows one to use a Newtons method to find the optimized values of $\sigma_A$ and $\Lambda$ with a few interactions. The uncertainties defined by Eq. (\ref{eq:Walphabetaalt}) are then found analytically. 

One should be careful in using smaller values of the convergence parameter, especially when there are a large number of parameters. The contribution to the variance of the emulator for a specific term in the polynomial expansion scales as $1/\Lambda^{2n}$, if $n$ is the order of the term. The number of terms of a given order scales as $N_{\rm pars}^n$. Thus, the importance of a given order relative to the previous order scales as $N_{\rm pars}/\Lambda^2$. If $\Lambda\lesssim N_{\rm pars}^{1/2}$, the emulator uncertainty does not converge well for increasing order. For example, for four model parameters, one might find that any model characterized by $\Lambda\gtrsim 2$ can be reasonably well emulated by a Taylor expansion, a model with 16 parameters might need $\Lambda\gtrsim 4$ to converge similarly as function of the order $K$. The challenge of handling emulation with large numbers of parameters is not confined to emulators based on polynomial expansions, but can be a problem for all emulators. A discussion of these issues is given in Sec. \ref{sec:optimize}.

\subsection{Contribution to the Emulator Uncertainty from the Uncertainty of the Hyper Parameters}
The emulator uncertainty $\sigma_Y(\vec{\theta})$ has contributions from the variance of the expansion coefficients for given values of $\sigma_A$ and $\Lambda$, and also contributions from the fact that $\sigma_A$ and $\Lambda$ are not known. Here, we present a general discussion of how one might approximate the latter contribution. We consider a set of hyper parameters, $\vec{h}$, which in this case would refer to $\sigma_A$ and $\Lambda$.

In this derivation, $\E(\mathcal{O}|\vec{h})$ refers to the the average of $\mathcal{O}$ over all the possible expansion coefficients given the training data for a specific set of hyper parameters $\vec{h}$. The notation $\E(\mathcal{O})$ refers to the average once one has also averaged over $\vec{h}$. Here, $\vec{h}_{0}$ is the hyper parameter vector at which the probability $P(\vec{h})$ is maximized, $\bar{Y}(\vec{h})=\E(Y|\vec{h})$, and $\bar{Y}_{0}=\E(Y|\vec{h}=\vec{h}_{0})$ and $\Var_0(Y)=\Var(Y|\vec{h}=\vec{h}_0)$. The goal is to estimate the variance of the emulated value at some point in parameter space, $\vec{\theta}$.
\begin{eqnarray}
\Var(Y) &=&
\int d\vec{h}~P(\vec{h}) \Var(Y|\vec{h})\\
\nonumber
&=&\int d\vec{h}~P(\vec{h}) \E([(Y-\bar{Y}(\vec{h})+(\bar{Y}(\vec{h})-\bar{Y}_{0})]^2|\vec{h})\\
\nonumber
&=&\int d\vec{h}~P(\vec{h})\left\{
\E((Y-\bar{Y}(\vec{h}))^2|\vec{h})+(\bar{Y}(\vec{h})-\bar{Y}_{0})^2)
\right\}\\
\nonumber
&=&\int d\vec{h}~P(\vec{h})\left\{\Var(Y|\vec{h}=\vec{h}_0)
+\frac{(h_i-h_{0i})(h_j-h_{0j})}{2}\frac{\partial^2}{\partial h_i\partial h_j}\Cov(Y|\vec{h})\right.\\
\nonumber
&&\hspace*{75pt}\left.+(h_i-h_{0i})(h_j-h_{0j})\left(\frac{\partial\bar{Y}(\vec{h})}{\partial h_i}\right)\left(\frac{\partial\bar{Y}(\vec{h})}{\partial h_j}\right)+\mathcal{O}(h-h_0)^3
\right\}\\
\label{eq:Sigma2FromW}
&\approx&\Var_0(Y)
+\left[\frac{\partial\bar{Y}(h)}{\partial h_i}\frac{\partial\bar{Y}(h)}{\partial h_j}\right]W_{ij},\\
\nonumber
W_{ij}&=&\int d\vec{h}~P(\vec{h})(h_i-h_{0i})(h_j-h_{0j}).
\end{eqnarray}
The approximation of the last step is built on the assumption that sufficient training data has been provided so that 
\begin{enumerate}\itemsep=0pt
\item $P(\vec{h})$ is sufficiently narrow about $\vec{h}_0$ that higher order terms in $\vec{h}-\vec{h}_0$ can be ignored.
\item $\Var(Y|\vec{h})$ is small so that terms that scale both as $\Var(Y|\vec{h})$ and $(\vec{h}-\vec{h}_0)^2$ can be ignored.
\end{enumerate}
One of the remaining terms is small because $\rm{Cov}(Y|\vec{h})$ is already small and additional powers of $\vec{h}-\vec{h}_0$ make it increasingly smaller compared to the first term. The neglected terms should all be increasingly smaller relative to the included terms as the number of training points is increased.

Applying the approximation of Eq. (\ref{eq:Sigma2FromW}) to the calculations of the uncertainty of the emulator,
\begin{eqnarray}\label{eq:sigmaEtot}
\sigma^2_Y(\vec{\theta})&\approx&\sigma_A^2\left[
C_0(\vec{\theta},\vec{\theta})
- C_0(\vec{\theta},\vec{\theta}_a)B_{ab}^{-1}C_0(\vec{\theta}_b,\vec{\theta})
\right]\\
\nonumber
&+&W_{\sigma_A\sigma_A}\left(\frac{\partial Y}{\partial{\sigma_A}}\right)^2
+2W_{\sigma_A\Lambda}\left(\frac{\partial Y}{\partial{\sigma_A}}\right)\left(\frac{\partial Y}{\partial{\Lambda}}\right)
+W_{\Lambda\Lambda}\left(\frac{\partial Y}{\partial{\Lambda}}\right)^2.
\end{eqnarray}

The r.h.s. of Eq. (\ref{eq:sigmaEtot}) is evaluated with $\sigma_A$ and $\Lambda$ chosen to maximize $Z(\sigma_a,\Lambda)$, and the covariance of the hyper-parameters, $W$, is approximated from the curvature of $Z(\sigma_A^2,\Lambda)$, as described in Eq. (\ref{eq:Walphabetaalt}). The first term represents the emulator's uncertainty due to the uncertainty of the coefficients $A_{\vec{n}}$ for fixed $\sigma_A$ and $\Lambda$, and the terms proportional to $W$ provide the uncertainty due to the fact that $\sigma_A$ and $\Lambda$ are themselves uncertain. In the limit of large numbers of training points the range of $\Lambda$ is well constrained around the optimum value and this becomes an excellent approximation, but the approximation becomes dubious when $\Lambda$ is poorly constrained. Often, the contributions to the uncertainty from the terms proportional to $W_{\alpha\beta}$ are ignored. In our studies it was found these corrections were relatively small when one chooses a reasonable set of training-point locations. However, these terms were found to be critical for optimizing the choice of locations of training points, $\vec{\theta}_{a=1\cdots N_{\rm train}}$. For example, if all the training points were placed very close to the origin, the terms were found to be significant, but for optimized placement the terms were very small. In contrast the strength of the other terms was observed to be only weakly dependent of the training-point placement. Thus, these additional terms played an outsized role when optimizing training-point placement.

Restating Eq. (\ref{eq:exactydely2}), the emulator depends on the training values as,
\begin{eqnarray}
Y(\vec{\theta})&=&C_0(\vec{\theta},\vec{\theta}_a)B^{-1}_{ab}y_b,
\end{eqnarray}
and is independent of $\sigma_A$ because $C_0$ and $B$ depend only on $\Lambda$. Thus, $\partial_{\sigma_A}Y(\vec{\theta})=0$, and one need only worry about $W_{\Lambda\Lambda}$ and can ignore $W_{\sigma_A\sigma_A}$ or $W_{\sigma_A,\Lambda}$, even though one needs to calculate all the terms in $W^{-1}$ before finding $W$. 

The probability $Z$ depends on $\sigma_A$, $\Lambda$, the location of the training points and the training values.
\begin{eqnarray}
\ln(Z)&=&-N_{\rm train}\ln(\sigma_A)-\frac{1}{2}\ln(|B|)-\frac{1}{2\sigma_A^2}y_aB^{-1}_{ab}y_b.
\end{eqnarray}
Because $B$ is independent of $\sigma_A$,
\begin{eqnarray}\label{eq:Winv0}
W^{-1}_{\sigma_A\sigma_A}&=&\frac{-N_{\rm train}}{\sigma_A^2}+\frac{3}{\sigma_A^4}y_aB^{-1}_{ab}y_b,\\
\nonumber
W^{-1}_{\sigma_A\Lambda}&\approx&\frac{1}{\sigma_A^3}y_a\left(\frac{d}{d\Lambda}B^{-1}\right)_{ab}y_b,\\
\nonumber
W^{-1}_{\Lambda\Lambda}&\approx&\frac{1}{2}\frac{d^2}{d\Lambda^2}\ln(|B|)-\frac{1}{2\sigma_A^2}y_a\left(
\frac{d}{d\Lambda}\right)^2B^{-1}_{ab}y_b.
\end{eqnarray}
These expressions assume that near the optimum values, $Z(\sigma_A,\Lambda)$ can be treated as if it were a Gaussian, so that the widths are determined by the curvature of $Z(\sigma_A,\Lambda)$ around the optimum. To take the derivatives w.r.t. $B^{-1}$, one can use the identity
\begin{eqnarray}
\frac{d}{d\Lambda}B^{-1}_{ab}&=&-B^{-1}_{ac}\left(\frac{d}{d\Lambda}B_{cd}\right)B^{-1}_{db}.
\end{eqnarray}
After defining
\begin{eqnarray}\label{eq:Bprimedef}
B'_{ab}&=&\Lambda^3\frac{d}{d\Lambda}B_{ab}=(|\vec{\theta}_a-\vec{\theta}_b|^2)B_{ab},\\
\nonumber
B''_{ab}&=&\Lambda^3\frac{d}{d\Lambda}B'_{ab}=(|\vec{\theta}_a-\vec{\theta}_b|^4)B_{ab}.
\end{eqnarray}
one finds
\begin{eqnarray}\label{eq:Winv}
W^{-1}_{\sigma_A\sigma_A}&=&N_{\rm train}\frac{1}{\sigma_A^2}-\frac{3}{\sigma_A^4}y_aB^{-1}_{ab}y_b,\\
\nonumber
W^{-1}_{\sigma_A\Lambda}&=&\frac{-1}{\Lambda^3\sigma_A^3}y_a\left(B^{-1}B'B^{-1}\right)_{ab}y_b,\\
\nonumber
W^{-1}_{\Lambda\Lambda}&=&\frac{1}{2}\frac{d^2}{d\Lambda^2}\ln(|B|)
+\frac{1}{2\sigma_A^2\Lambda^6}y_a\left(
2B^{-1}B'B^{-1}B'B^{-1}-B^{-1}B''B^{-1}
\right)_{ab}y_b\\
\nonumber
&+&\frac{3}{2\sigma_A^2\Lambda^4}y_a\left(B^{-1}B'B^{-1}\right)_{ab}y_b.
\end{eqnarray}
To calculate $(d^2/d\Lambda^2)\ln(|B|)$, one can use the identity
\begin{eqnarray}
\frac{d}{d\Lambda}|B|&=&|B|{\rm tr}\left(B^{-1}\frac{dB}{d\Lambda}\right).
\end{eqnarray}
This leads to
\begin{eqnarray}
  \frac{d^2}{d\Lambda^2}\ln(|B|)&=&
  -{\rm tr}\left(B^{-1}\frac{dB}{d\Lambda}B^{-1}\frac{dB}{d\Lambda}\right)
  +{\rm tr}\left(B^{-1}\frac{d^2B}{d\Lambda^2}\right).
\end{eqnarray}

The last piece required to calculate the emulator uncertainty in Eq. (\ref{eq:sigmaEtot}) is $\partial_\Lambda Y(\vec{\theta})$.
\begin{eqnarray}\label{eq:dEdLambda}
\partial_\Lambda Y(\vec{\theta})&=&
\partial_\Lambda C_0(\vec{\theta},\vec{\theta}_a)B^{-1}_{ab}y_b\\
\nonumber
&=&
\left\{\frac{1}{\Lambda^3}|\vec{\theta}-\vec{\theta}_a|^2C_0(\vec{\theta},\vec{\theta}_a)B^{-1}_{ab}
-C_0(\vec{\theta},\vec{\theta}_a)(B^{-1}\left[\frac{d}{d\Lambda}B\right]B^{-1})_{ab}
\right\}y_b.
\end{eqnarray}
Thus, to estimate the contribution to the emulator uncertainty due to the uncertainties in determining the hyper parameters one should perform the following steps.
\begin{enumerate}\itemsep=0pt
\item First tune the emulator and find the optimum values of $\sigma_A$ and $\Lambda$.
\item Find the $2\times 2$ matrix $W^{-1}$ by using Eq.s (\ref{eq:Winv}) with the matrices $B,B'$ and $B''$ calculated using $\Lambda$ from step 1. This also involves finding $(d^2/d\Lambda^2)\ln(|B|)$.  
\item Invert $W^{-1}$ to find $W$. Only $W_{\Lambda\Lambda}$ will be needed because $\partial Y/\partial\sigma_A=0$.
\item Find $dY/d\Lambda$ from Eq. (\ref{eq:dEdLambda}) and using $W_{\Lambda\Lambda}$ calculate the contribution to $\sigma_Y^2$ in Eq. (\ref{eq:sigmaEtot}).
\end{enumerate}

\subsection{The Emulator Uncertainty at Training Points}

First, we demonstrate the well-known result that in the absence of point-by-point noise, i.e. $\Delta=0$, the estimated uncertainty will zero. If one sets $\vec{\theta}=\vec{\theta}_c$, one of the training points, and if $\Delta=0$, then
\begin{eqnarray}
C_0(\vec{\theta}_a,\vec{\theta}_c)&=& B_{ac},\\
\nonumber
T_{\vec{m}}(\vec{\theta}_c)T_{\vec{m}}(\vec{\theta}_c)&=& B_{cc},~{\rm no~sum~over~}c,
\end{eqnarray}
and
\begin{eqnarray}
Y(\vec{\theta}_c)&=&\chi_a
C_0(\vec{\theta}_a,\vec{\theta}_c)=B^{-1}_{ab}y_b B_{ac}=y_c,\\
\nonumber
\sigma^2_Y(\vec{\theta}_c)|_{\Delta=0}
&=&\sigma_A^2 B_{cc}-\sigma_A^2 B_{ca}B^{-1}_{ab} B_{bc}+W_{\Lambda\Lambda}\left(\frac{\partial Y}{\partial\Lambda}\right)^2\\
\nonumber
&=&W_{\Lambda\Lambda}\left(\frac{\partial Y}{\partial\Lambda}\right)^2.
\end{eqnarray}

Next, using Eq. (\ref{eq:dEdLambda}), one can see that evaluated at a training point, $\vec{\theta}_c$, is
\begin{eqnarray}
\frac{\partial Y(\vec{\theta}_c)}{\partial\Lambda}&=&\left[B'_{ca}B^{-1}_{ab}-B_{ca}(B^{-1}B'B^{-1})_{ab}\right]y_b=0.
\end{eqnarray}
Here, $B'=(d/d\Lambda)B$. Thus, if $\Delta=0$, there is no uncertainty when evaluated at the training points and the optimum value for the emulator exactly reproduces the training values. 

When $\Delta\ne 0$, the optimum is no longer constrained to exactly reproduce the full model at the training points, and the uncertainty at a training point, $\vec{\theta}_c$ is no longer zero,
\begin{eqnarray}\label{eq:emulatorerror}
{\rm Var}(Y(\vec{\theta}_c))&=&
\sigma_A^2\left(\Delta_{cc}-\Delta_{ca}B^{-1}_{ab}\Delta_{bc}\right),~~~({\rm no~sum~over~} c). 
\end{eqnarray}
If $\Delta$ is small, the uncertainty of the emulator is simply given by $\sigma_A \Delta^{1/2}$, however it is somewhat reduced for larger $\Delta$. 

Typically, the emulator will be applied thousands, or perhaps millions of times, all for different $\vec{\theta}$. Before evaluating the emulator, it would be efficient to first calculate and store both $B^{-1}_{ab}$ and $\chi_{a}$. Calculating the matrix $B$ involves calculating $N_{\rm train}(N_{\rm train}+1)/2$ exponentials, which is extremely fast. One must then invert $B$, an $N_{\rm train}\times N_{\rm train}$ matrix. Then, when considering a specific $\vec{\theta}$, one calculates $C_0(\vec{\theta}_a,\vec{\theta})$ for each $a$ according to Eq. (\ref{eq:C0default}).  Once these are stored, calculating $\rm{Var}(Y(\vec{\theta}))$ involves calculating a double sum with $N_{\rm train}^2$ elements. 

The derivations above could have been performed finding the coefficients $A_{\vec{n}}$ that maximize the probability, $Z(\vec{A})$, rather than solving for the mean coefficients, but the expressions for $Y(\vec{\theta})$ and $\sigma_Y(\vec{\theta})$ end up identical to those provided above.

\subsection{Summarizing the Emulation}

Despite the fact that the emulator was built in terms of optimized expansion coefficients $A_{\vec{n}}$, the entire emulator construction and application avoids any mention of the coefficients. Through the kernel trick, the emulator need only know the form for the unconstrained correlation, $C_0(\vec{\theta}_1,\vec{\theta}_2)$. Using the short hand definitions of the $N_{\rm train}\times N_{\rm train}$ matrix and $N_{\rm train}$-dimensional vector $\chi_a$,
\begin{eqnarray}
B_{ab}&\equiv&C_0(\vec{\theta}_a,\vec{\theta}_b)+\Delta_{ab},\\
\nonumber
\chi_a&\equiv&B^{-1}_{ab}y_b,
\end{eqnarray}
one finds expressions for the following quantities,
\begin{eqnarray}
\label{eq:summary}
\sigma_A^2&=&\frac{1}{N_{\rm train}}\chi_ay_a,\\
\nonumber
Z(\sigma_A^2,\Lambda)&=&\sigma_A^{-N_{\rm train}}
|B|^{-1/2}\exp\left[
-\frac{1}{2\sigma_A^2}y_aB^{-1}_{ab}y_b\right]\\
\nonumber
&=&\sigma_A^{-N_{\rm train}}|B|^{-1/2}e^{-1/(2N_{\rm train})}~~{\rm when~}\sigma_a~{\rm is~optimized},\\
\nonumber
Y(\vec{\theta})&=&\chi_aC_0(\vec{\theta}_a,\vec{\theta}),\\
\nonumber
\sigma^2_Y(\vec{\theta})&=&\sigma_A^2\left[
C_0(\vec{\theta},\vec{\theta})
- C_0(\vec{\theta},\vec{\theta}_a)B_{ab}^{-1}C_0(\vec{\theta}_b,\vec{\theta})
\right]+W_{\Lambda\Lambda}\left(\frac{\partial Y}{\partial{\Lambda}}\right)^2.
\end{eqnarray}
Aside from the last term, proportional to $W_{\Lambda\Lambda}$, the other pieces are no different than what finds with GP emulators with Gaussian correlation kernels.

Outlining the procedure for tuning,
\begin{enumerate}\itemsep=0pt
\item Choose the locations of the $N_{\rm train}$ training points, $\vec{\theta}_a$.
\item Choose a value of $\Lambda$.
\item Evaluate the full model at the training points, i.e. determine $y_{a=1,N_{\rm train}}$.
\item Define the noise matrix $\Delta_{ab}$. For example, if the noise is uncorrelated, and is approximately one percent of the signal, one would choose $\Delta_{ab}=\alpha^2\delta_{ab}$, with $\alpha=0.01$.
\item Calculate and store $B$, $B_{\rm inv}$, $\chi$.
\item Calculate $\sigma_A^2$ according to Eq. (\ref{eq:sigmaAguess}).
\item Calculate $Z(\sigma_A^2,\Lambda)$.
\item Repeat steps 2-7 for different $\Lambda$ and keep the value of $\Lambda$, and its corresponding $\sigma_A$,  that maximized $Z$.
\item Once both $\Lambda$ and $\sigma_A$ are determined, calculate and store $W_{\Lambda\Lambda}$.
\end{enumerate}
To this point, the most numerically taxing part of the procedure is inverting $B$ and calculating its determinant. One might need to repeat the steps, as mentioned in (8) roughly a dozen times to best find $\Lambda$.

Once this is finished, one can evaluate the emulator and its uncertainty at any $\vec{\theta}$. At this point one needs to merely calculate a vector $C_0(\vec{\theta}_a,\vec{\theta})$ for $a=1-N_{\rm train}$, then contract the vector with $\chi$ to obtain $Y(\vec{\theta})$ or contract with $B^{-1}$ to find the squared emulator uncertainty.

\setcounter{section}{2}

\section{Choosing the Optimum Training Points}
\label{sec:optimize}

For emulating a convergent function, optimizing one's choice of training points involves choosing points that allow accurate determination of the lowest-order expansion coefficients. If one has $N_{\rm pars}$ model parameters, one needs at least $N_{\rm train}=N_{\rm pars}+1$ points to fit all the up-to-order-one coefficients, i.e. the slope in $N_{\rm pars}$ dimensions and the intercept. To obtain the coefficients up to second order, one would need at least $(N_{\rm pars}+1)(N_{\rm pars}+2)/2$ training values. Thus, it might seem natural to choose a set of either $N_{\rm pars}+1$ or $(N_{\rm pars}+1)(N_{\rm pars}+2)/2$ training points with optimized positions. 

One would like to find an algorithm than optimizes the constraining power of the training points for arbitrary priors, for an arbitrary numbers of model parameters, and for an arbitrary number of training points. This needs to be performed before the full model is run, and thus before any of the training values, $y_a$, are known. This means that the hyper-parameters, $\Lambda$ and $\sigma_A$, are also as yet undetermined. In this section, algorithms are presented that accomplish this goal, but as will be shown, depend on an estimate for $\Lambda$. Given $\Lambda$, an optimized placement of the training points is found. Fortunately, the placement will not be strongly dependent on $\Lambda$. Unfortunately, the overall uncertainties are strongly sensitive to $\Lambda$, so choosing the number of training points required for a specific accuracy can be problematic.

An accuracy metric is defined based on the expected squared uncertainty averaged over the prior. If the priors are all some combination of Gaussian or uniform forms, analytic expressions are found for this metric. Although these expressions are rather messy, they enable a rather rapid numerically optimized positioning of the training points. The next subsection presents the analytic expressions for the accuracy metric. The analytic forms require that the prior can be expressed in terms of a product of either uniform or Gaussian priors for each model parameter. The next subsection provides a case study of how the accuracy behaves with the number of training points, with the number of parameters, and as a function of the estimated convergence parameter. For context, results are compared to Latin hyper-cube sampling \cite{lhs}. 

Certain numbers of training points allow the placement to satisfy symmetries of the prior, if such symmetries exist. If the prior is uniform within a $N_{\rm train}-$dimensional cube, the training points might be placed at the corners of an interior cube if $N_{\rm train}=2^{N_{\rm pars}}$, or perhaps adding an additional point at the center. For spherically symmetric priors with $N_{\rm train}=N_{\rm pars}+1$, one might place the points in a simplex. The additional improvement for these specific values is rather modest, and is discussed in Sec. \ref{subsec:finestructure}.

\subsection{Expressing the Overall Emulator Accuracy}
 
To optimize the location of the training points, we define the accuracy metric as  being proportional to the emulator's uncertainty squared, $\sigma^2_Y(\vec{\theta})$, averaged over the prior, $P(\vec{\theta})$, then divided by $\sigma_A^2$. One might, or might not, be able to make solid arguments that averaging $\sigma_Y$ to some other power might be more appropriate. This choice has the benefit of allowing much of the averaging to be performed analytically, and seems to provide a reasonable balance of weighing areas with larger uncertainty relative to areas with lower uncertainty.

The goal is thus to calculate the expected value of $\sigma_A^2(\vec{\theta})$, then average this weighted by the prior, $P(\vec{\theta})$. This averaging over $\vec{\theta}$ will be denoted with the notation, $\mathbf{E}_P(\sigma_Y^2)$. Because the hyper-parameters are not yet determined at the time when the training points are chosen, the goal will be to minimize $\mathbf{E}_P(\sigma_Y^2)$ for some assumed values of $\sigma_A$ and $\Lambda$.
\begin{eqnarray}\label{eq:sigmabar}
\mathbf{E}_P(\sigma_Y^2)&\equiv&\int d\vec{\theta} P(\vec{\theta}) \sigma^2_Y(\vec{\theta})\\
\nonumber
&=&\int d\vec{\theta} P(\vec{\theta})\left\{
\sigma_A^2\left[1-C_0(\vec{\theta},\vec{\theta}_a)B^{-1}_{ab}C_0(\vec{\theta_b},\vec{\theta})\right]
+\left(\frac{\partial Y}{\partial\Lambda}\right)^2W_{\Lambda\Lambda}
\right\}.
\end{eqnarray}
The expressions from Sec. (\ref{sec:theory}) for $W_{\Lambda\Lambda}$ and $(\partial Y/\partial\Lambda)^2$ in Eqs. (\ref{eq:Winv}) and Eq. (\ref{eq:dEdLambda}) depend on the training values through the product $y_ay_b$. Given that the positions of the training points are not yet known, the quantities $y_ay_b/\sigma_a^2$ must be replaced with their unconstrained expectation,
\begin{eqnarray}\label{eq:yaybbar}
\frac{y_ay_b}{\sigma_A^2} \rightarrow B_{ab}&=&\left[C_0(\vec{\theta}_a,\vec{\theta}_b)+\alpha^2\delta_{ab}\right].
\end{eqnarray}
The expressions for $W$, which are independent of $\vec{\theta}$, then become
\begin{eqnarray}
W^{-1}_{\sigma_A\sigma_A}&=&\frac{2N_{\rm train}}{\sigma_A^2},\\
\nonumber
W^{-1}_{\sigma_A\Lambda}&=&\frac{-1}{\Lambda^3\sigma_A}{\rm Tr}(B'B^{-1}),\\
\nonumber
W^{-1}_{\Lambda\Lambda}&=&
+\frac{1}{2\Lambda^6}{\rm Tr}(B'B^{-1}B'B^{-1}).
\end{eqnarray}
One needs to invert $W^{-1}$ to find $W_{\Lambda\Lambda}$, which will be independent of $\sigma_A$.

The factor $(\partial_\Lambda Y)^2$ depends on $\vec{\theta}$ and becomes
\begin{eqnarray}
\label{eq:d2EdL2}
(\partial_\Lambda Y)^2(\vec{\theta})&=&\frac{\sigma_A^2}{\Lambda^6}\sum_{ab}\left\{
C_0(\vec{\theta}_b,\vec{\theta})|\vec{\theta}_b-\vec{\theta}|^2|\vec{\theta}-\vec{\theta}_a|^2
C_0(\vec{\theta},\vec{\theta}_a)B^{-1}_{ab}
\right\}\\
\nonumber
&-&\frac{\sigma_A^2}{\Lambda^6}\sum_{ab}\left\{
C_0(\vec{\theta}_b,\vec{\theta})|\vec{\theta}-\vec{\theta}_a|^2C_0(\vec{\theta},\vec{\theta}_a)(B^{-1}B'B^{-1})_{ab}
\right\}\\
\nonumber
&-&\frac{\sigma_A^2}{\Lambda^6}\sum_{ab}\left\{
C_0(\vec{\theta}_b,\vec{\theta})|\vec{\theta}-\vec{\theta}_b|^2C_0(\vec{\theta},\vec{\theta}_a)(B^{-1}B'B^{-1})_{ab}
\right\}\\
\nonumber
&+&\frac{\sigma_A^2}{\Lambda^6}\sum_{ab}\left\{
C_0(\vec{\theta}_b,\vec{\theta})C_0(\vec{\theta},\vec{\theta}_a)(B^{-1}B'B^{-1}B'B^{-1})_{ab}
\right\}.
\end{eqnarray}
The expression for $\sigma_Y^2$ above has a simple dependence w.r.t $\sigma_A$. Every term is proportional to $\sigma_A^2$. In contrast, the  $\Lambda$ dependence is complicated. To optimize the positions of the training points, $\vec{\theta}_{a=1\cdots N_{\rm train}}$, one must minimize the r.h.s. of Eq. (\ref{eq:sigmabar}). Unfortunately, without any input from full-model runs, one can only guess at $\Lambda$. Before training all values of $\sigma_A^2$ are considered equally likely. Thus, because of the simple proportionality to $\sigma_A^2$, one can choose any $\sigma_A$ and it will not affect the choice of training-point locations to minimize $\mathbf{E}_P(\sigma_Y^2)$. Thus, as a metric for the expected accuracy, we will choose
\begin{eqnarray}\label{eq:accuracymetric}
{\rm Accuracy~metric}&=&\frac{\mathbf{E}_P(\sigma_Y^2)}{\sigma_A^2},
\end{eqnarray}
which is independent of $\sigma_A$ and depends only on $\Lambda$ and on the location of the training points. The goal will be to find locations of the training points that minimize the metric after guessing at a value for $\Lambda$.

Once $\Lambda$ is chosen, one can integrate Eq. (\ref{eq:sigmabar}) over the prior to calculate the accuracy. Because the matrices, $B$, $B'$ and $B''$, defined in Eq.s (\ref{eq:Bdef}) and (\ref{eq:Bprimedef}), depend only on $\Lambda$ and on the locations of the training points, one can express the metric in terms of the integrals, $I_{ab}, J_{ab}$ and $K_{ab}$ defined as:
\begin{eqnarray}\label{eq:sigma2prior}
\frac{\mathbf{E}_P(\sigma_Y^2)}{\sigma_A^2}&=&{\rm Tr}\left(\unity-IB^{-1}\right)+\frac{W_{\Lambda\Lambda}}{\Lambda^6}
{\rm Tr}\left\{
KB^{-1}+JB^{-1}B'B^{-1}+IB^{-1}B'B^{-1}B'B^{-1}
\right\},\\
\nonumber
I_{ab}&\equiv&\int d\vec{\theta}~P(\vec{\theta})C_0(\vec{\theta}_a,\vec{\theta})C_0(\vec{\theta},\vec{\theta}_b),\\
\nonumber
J_{ab}&\equiv&J_{ab;a}+J_{ab;b},\\
\nonumber
J_{ab;a}&\equiv&\int d\vec{\theta}~P(\vec{\theta})C_0(\vec{\theta}_a,\vec{\theta})
|\vec{\theta}_a-\vec{\theta}|^2C_0(\vec{\theta},\vec{\theta}_b),\\
\nonumber
K_{ab}&\equiv&\int d\vec{\theta}~P(\vec{\theta})C_0(\vec{\theta}_a,\vec{\theta})
|\vec{\theta}_a-\vec{\theta}|^2|\vec{\theta}-\vec{\theta}_b|^2C_0(\vec{\theta},\vec{\theta}_b).
\end{eqnarray}
Here, $a$ and $b$ run from 1 to $N_{\rm train}$ The $N_{\rm train}\times N_{\rm train}$ matrices $I,J$ and $K$ carry all the relevant information about the prior.

If the prior $P(\vec{\theta})$ can be written as product of priors for individual model parameters,
\begin{eqnarray}
P(\vec{\theta})&=&\prod_i P_i(\theta_i),
\end{eqnarray}
one can express the integrals over the prior in terms of integrals over individual parameters.
\begin{eqnarray}
I_{ab}&=&\prod_i I_{ab;i},\\
\nonumber
J_{ab}&=&I_{ab}\sum_i\frac{J_{ab;ai}+J_{ab;bi}}{I_{ab;i}},\\
\nonumber
K_{ab}&=&\frac{J_{ab;a}J_{ab;b}}{I_{ab}}+I_{ab}\sum_i\left\{
\frac{K_{ab;i}}{I_{ab;i}}-\frac{J_{ab;ai}J_{ab;bi}}{I_{ab;i}^2}
\right\}.
\nonumber
\end{eqnarray}
Here, the integrals over single parameters are:
\begin{eqnarray}\label{eq:IJKi}
I_{ab;i}&=&\int d\theta_i~P_i(\theta_i)\exp\left\{-\frac{(\theta_i-\theta_{ai})^2}{2\Lambda^2}-\frac{(\theta_i-\theta_{bi})^2}{2\Lambda^2}
\right\},\\
\nonumber
J_{ab;ai}&=&\int d\theta_i~P_i(\theta_i)(\theta_i-\theta_{ai})^2\exp\left\{-\frac{(\theta_i-\theta_{ai})^2}{2\Lambda^2}-\frac{(\theta_i-\theta_{bi})^2}{2\Lambda^2}
\right\},\\
\nonumber
K_{ab;i}&=&\int d\theta_i~P_i(\theta_i)(\theta_i-\theta_{ai})^2(\theta_i-\theta_{bi})^2\exp\left\{-\frac{(\theta_i-\theta_{ai})^2}{2\Lambda^2}-\frac{(\theta_i-\theta_{bi})^2}{2\Lambda^2}
\right\}.
\end{eqnarray}
After calculating the integrals $I_{ab;i}$ analytically, the integrals $J_{ab;ai}, J_{ab;bi}$ and $K_{ab;i}$ can all be generated by taking derivatives of $I_{ab;i}$ w.r.t. $\gamma=1/\Lambda^2$.

One also needs an expression for $W$. This is found in the same way as for Eq.s (\ref{eq:Winv}). In those expressions $W_{\rm inv}$ was defined in terms of $B_{ab}$, where $B_{ab}$ was expressed in terms of the training values through the matrix $y_ay_b$, whereas in this case $B_{ab}$ is estimated from the pre-constrained covariance, defined in Eq. (\ref{eq:yaybbar}).

\subsubsection{Analytic expressions for Gaussian and uniform priors}
If a specific model parameter, $\theta_i$, has a Gaussian prior,
\begin{eqnarray}\label{eq:guassianbeta}
P_i(\theta=\theta_i)&=&\frac{1}{(2\pi/\beta_i)^{1/2}}e^{-\beta_i\theta^2/2},
\end{eqnarray}
one can integrate Eq.s (\ref{eq:IJKi}) and find
\begin{eqnarray}
I_{ab;i}&=&\sqrt{\frac{\beta_i}{\lambda_i}}\exp
\left\{
-\frac{1}{2\lambda_i}[\gamma^2(\theta_a-\theta_b)^2+\beta_i\gamma(\theta_a^2+\theta_b^2)]
\right\},\\
\nonumber
J_{ab;ai}&=&I_{ab;i}\left\{
\frac{1}{\lambda_i}+|\vec{\theta}_a-\vec{\theta}_b|^2\frac{\gamma^2+\beta_i\gamma}{\lambda_i^2}
+(\theta_a^2+\theta_b^2)\frac{\beta_i^2}{2\lambda_i^2}
+\frac{\beta_i}{2\lambda_i}(\theta_a^2-\theta_b^2)
\right\},\\
\nonumber
J_{ab;bi}&=&I_{ab;i}\left\{
\frac{1}{\lambda_i}+|\vec{\theta}_a-\vec{\theta}_b|^2\frac{\gamma^2+\beta_i\gamma}{\lambda_i^2}
+(\theta_a^2+\theta_b^2)\frac{\beta_i^2}{2\lambda_i^2}
-\frac{\beta_i}{2\lambda_i}(\theta_a^2-\theta_b^2)
\right\},\\
\nonumber
K_{ab;i}&=&\frac{J_{ab;ai}J_{ab;bi}}{I_{ab;i}}+I_{ab;i}\left\{
\frac{2}{\lambda_i^2}
+2\beta_i^2(\theta_a^2+\theta_b^2)/\lambda_i^3
-(4\gamma^2+2\beta_i\lambda_i)(\theta_a-\theta_b)^2/\lambda_i^3
\right\},
\end{eqnarray}
where $\lambda_i\equiv2\gamma+\beta_i$ and $\gamma=1/\Lambda^2$.

Similarly, one can perform the integrals for uniform priors. In this case the prior is
\begin{eqnarray}\label{eq:uniformthetamax}
P_i(\theta=\theta_i)&=&\left\{\begin{array}{rl}
\frac{1}{2\theta_{{\rm max},i}},&-\theta_{{\rm max},i}<\theta<\theta_{{\rm max},i}\\
0,&{\rm otherwise}\\
\end{array}\right.
\end{eqnarray}
After making the definitions
\begin{eqnarray}
\bar{\theta}&\equiv&\frac{\theta_a+\theta_b}{2},\\
\nonumber
\Delta\theta&\equiv&\frac{\theta_a-\theta_b}{2},\\
\nonumber
\theta_+&\equiv&\theta_{{\rm max},i}-\bar{\theta},\\
\nonumber
\theta_-&\equiv&-\theta_{{\rm max},i}-\bar{\theta},\\
p&\equiv&\frac{1}{2\theta_{{\rm max},i}}e^{-\gamma\Delta\theta^2},\\
\nonumber
q&\equiv&\frac{1}{2}\sqrt{\pi}{\gamma}\left[
{\rm erf}(\sqrt{\gamma}\theta_+) - {\rm erf}(\sqrt{\gamma}\theta_-)
\right]\\
\nonumber
r&\equiv&\theta_+ e^{-\gamma\theta_+^2}-\theta_-e^{-\gamma\theta_-^2},
\end{eqnarray}
the integrals of interest are
\begin{eqnarray}
I_{ab;i}&=&pq,\\
\nonumber
J_{ab;ai}&=&(\frac{1}{2\gamma}-\Delta\theta^2)I_{ab;i}-\frac{1}{2\gamma}pr
+pe^{-\gamma\theta_+^2}\frac{\Delta\theta}{\gamma}-Pe^{-\gamma\theta_-^2}\frac{\Delta\theta}{\gamma},\\
\nonumber
J_{ab;bi}&=&(\frac{1}{2\gamma}-\Delta\theta^2)I_{ab;i}-\frac{1}{2\gamma}pr
-pe^{-\gamma\theta_+^2}\frac{\Delta\theta}{\gamma}+Pe^{-\gamma\theta_-^2}\frac{\Delta\theta}{\gamma},\\
\nonumber
K_{ab;i}&=&\frac{1}{2}(J_{ab;ai}+J_{ab;bi})
\left(\frac{1}{2\gamma}+\Delta\theta^2\right)\\
\nonumber
&+&I_{ab;i}\left(
\frac{1}{2\gamma^2}-\frac{2\Delta\theta^2}{\gamma}
\right)\\
\nonumber
&-&pe^{-\gamma\theta_+^2}\theta_+\left(\frac{1}{2\gamma^2}-\frac{3\Delta\theta^2}{2\gamma}+\frac{\theta_+^2}{2\gamma}\right)\\
\nonumber
&+&pe^{-\gamma\theta_-^2}\theta_-\left(\frac{1}{2\gamma^2}-\frac{3\Delta\theta^2}{2\gamma}+\frac{\theta_-^2}{2\gamma}\right).
\end{eqnarray}
For most of the results presented in this section Gaussian priors were considered. 

\subsection{Dependence on the number of parameters, the number of training points and the convergence parameter}

The accuracy metric, defined in Eq. (\ref{eq:accuracymetric}), was optimized by performing a simple search through the parameter space. For some initial configuration of training points $\mathbf{E}_P(\sigma_Y^2)/\sigma_A^2$ was calculated. A new set of training points was then chosen by taking a small random step in parameter space for each training points. If the new uncertainty was reduced, the step was accepted. This is by no means a particularly efficient method, but the method was sufficiently fast for $N_{\rm pars}\lesssim 10$, as the optimized values were found within a few minutes or less. Larger numbers of parameters required increasing amounts of patience.  

\begin{figure}
\begin{minipage}{0.475\textwidth}
\includegraphics[width=0.95\textwidth]{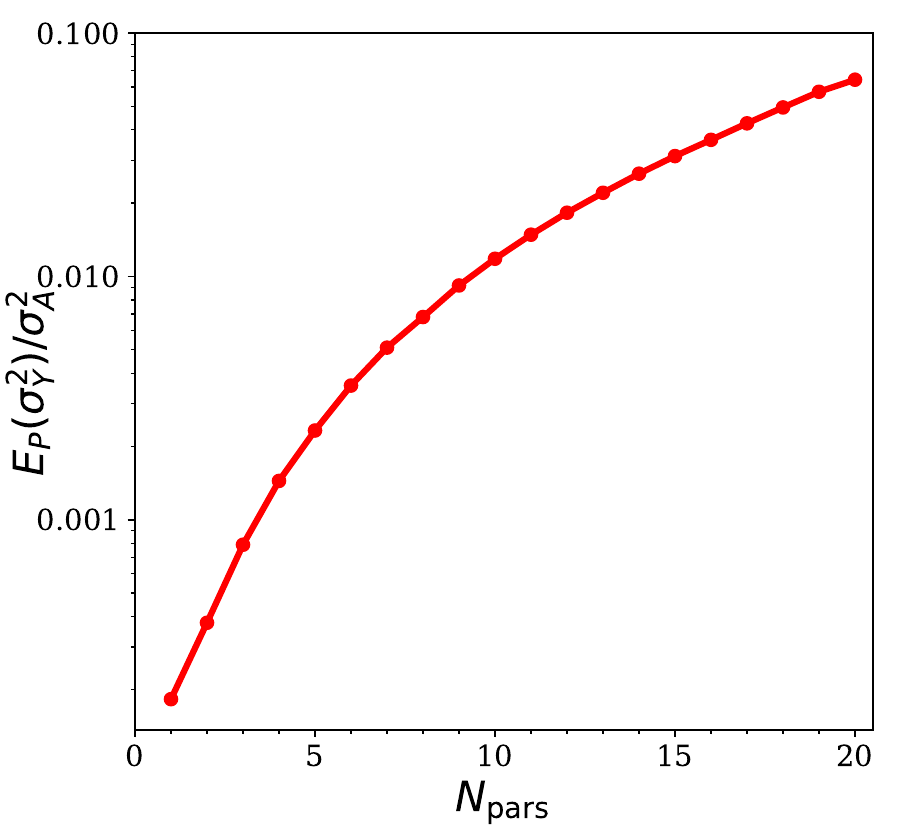}
\caption{\label{fig:sigma2vsNpars}
The accuracy metric as a function of the number of model parameters assuming Gaussian priors. The convergence parameter was set to $\Lambda=3$, the point-by-point noise was set to $\alpha=0.01$ and the number of training points was set to match the number of expansion coefficient of quadratic order or less. This illustrates the difficulty to accurately emulate higher dimensional functions even if the number of training points increases as $N_{\rm train}=(N_{\rm pars}+1)(N_{\rm pars}+2)/2$.}
\end{minipage}
\hspace*{0.025\textwidth}
\begin{minipage}{0.475\textwidth}
\includegraphics[width=\textwidth]{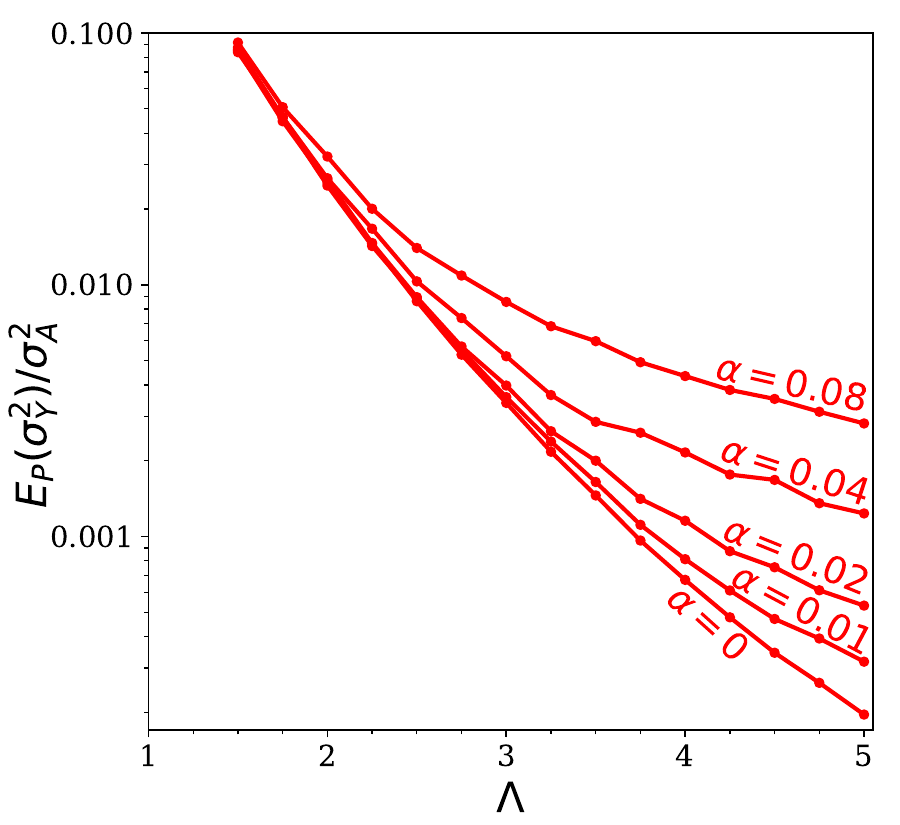}
\caption{\label{fig:Sigma2vsLambda}
The expected accuracy of the emulator falls rapidly with the convergence parameter $\Lambda$. For this illustration $N_{\rm pars}=6$ and the number of training points was set to 28, the number of expansion coefficients of order $n=2$ or less. Calculations were repeated for different values of the point-by-point noise, $\alpha$. As expected, the impact of the point-by-point noise sets in as the accuracy of the emulator improves.}
\end{minipage}
\end{figure}

The difficulty in emulating high-dimension parameter spaces increases more than linearly with the number of model parameters. This is because the number of expansion coefficients of order $K$ is
\begin{eqnarray}\label{eq:NcCounta}
N_{\rm coefficients}(K=0)&=&1,\\
\nonumber
N_{\rm coefficients}(K=1)&=&N_{\rm pars},\\
\nonumber
N_{\rm coefficients}(K=2)&=&N_{\rm pars}(N_{\rm pars}+1)/2,\\
\nonumber
N_{\rm coefficients}(K=3)&=&N_{\rm pars}(N_{\rm pars}+1)(N_{\rm pars}+2)/6,\\
\nonumber
\vdots
\end{eqnarray}
Further, the number of coefficients of order less than or equal to $K$, is
\begin{eqnarray}\label{eq:NcCountb}
N_{\rm coefficients}(K\le 0)&=&1,\\
\nonumber
N_{\rm coefficients}(K\le 1)&=&N_{\rm pars}+1,\\
\nonumber
N_{\rm coefficients}(K\le2)&=&(N_{\rm pars}+1)(N_{\rm pars}+2)/2,\\
\nonumber
N_{\rm coefficients}(K\le 3)&=&(N_{\rm pars}+1)((N_{\rm pars}+2)(N_{\rm pars}+3)/6,\\
\nonumber
\vdots
\end{eqnarray}
Thus, the number of terms of order $K+1$ is larger than the number of terms of order $K$ by a factor of $(N_{\rm pars}+K)/K$. If the number of parameters is large, there are large numbers of parameters of higher order. For example for $N_{\rm pars}=12$ there are 364 cubic terms, and for $N_{\rm pars}=20$ there are 1540 cubic terms. Because the importance of higher order terms falls as $1/\Lambda^{2n}$, emulating high-dimensional functions requires larger values of the convergence parameter to arrive at the same accuracy. Figure \ref{fig:sigma2vsNpars} shows the accuracy metric as a function of the number of parameters for fixed $\Lambda=3$ and for a number of training points set to correspond to the number of expansion parameters of order two or less, $N_{\rm train}=(N_{\rm pars}+1)(N_{\rm pars}+2)/2$.
 
It should be emphasized that this measure of the uncertainty is proportional to $\mathbf{E}_P(\sigma_Y^2)$, not to $\mathbf{E}_P(\sigma_Y)$. Thus, to get to 10\% accuracy, one might want $\mathbf{E}_P(\sigma_Y^2)/\sigma_A^2$ to be close to 0.01.

If $\Lambda$ is higher, it is easier to emulate functions because the contribution to $\mathbf{E}_P(\sigma_Y^2)$ of terms of order $n$ falls roughly as $1/\Lambda^2$. Figure \ref{fig:Sigma2vsLambda} shows how the expected accuracy falls as a function of $\Lambda$. For this illustration the number of parameters was set to 6 and the number of training points was again set to correspond to the number of expansion coefficients of order $n=2$ or less. Figure \ref{fig:Sigma2vsLambda} also illustrates the dependence on the point-by-point noise-strength, $\alpha$. The effect of $\alpha$ only sets in when the accuracy metric falls to near $\alpha^2$. However, given that there are multiple training points are considered the expected  uncertainty can fall below $\alpha^2$.

\begin{figure}
\begin{minipage}{0.475\textwidth}
\includegraphics[width=0.95\textwidth]{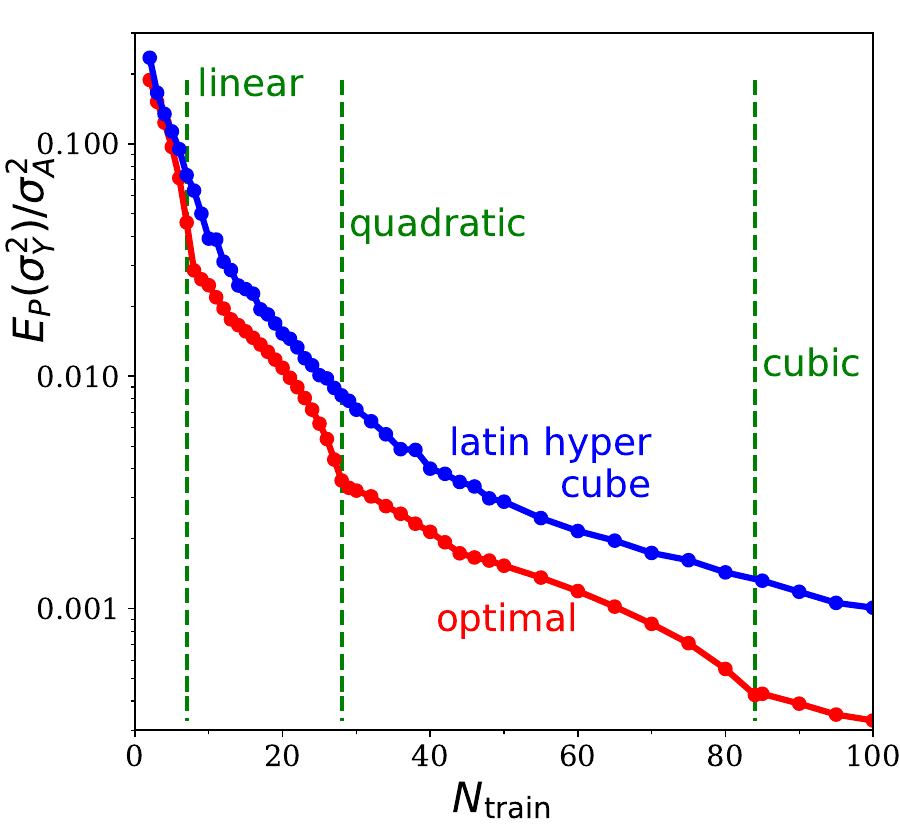}
\caption{\label{fig:Sigma2vsNTrain}
For all Gaussian priors, with $N_{\rm pars}=6$, $\Lambda=3$, and $\alpha=0.01$, the expected mean accuracy improves with increasing the number of training points. The minimum number of training points required for linear, quadratic and cubic fits are 7, 28 and 82. After passing one of these threhsolds the marginal improvement for additional training points decreases. For $N=8$, one can place the points in a simplex with an additional point at the center, and the symmetry of this configuration provides an especially large improvement. This behavior might inspire the modeler to choose a number of training points equal to one of the thresholds. The optimized placement provides significant improvement compared to Latin hyper-cube sampling.
}
\end{minipage}
\hspace*{0.025\textwidth}
\begin{minipage}{0.475\textwidth}
\includegraphics[width=\textwidth]{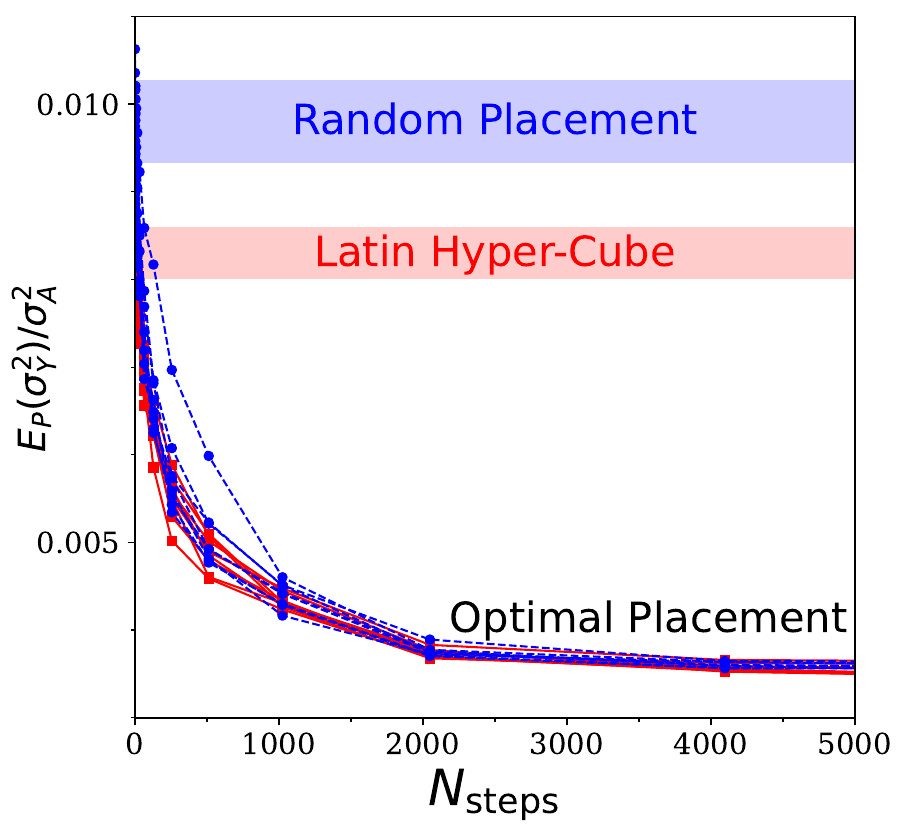}
\caption{\label{fig:Sigma2vsNMC}
The search procedure provides most of the improvement for the first few thousand steps. For this example, $N_{\rm pars}=6$, $\Lambda=3$ and $\alpha=0.01$ the number of training points was set to 28. Several initial configurations were considered. Those where the points were placed randomly (blue circles) tended to start at roughly three times the optimized value, whereas those that were initialized with Latin hyper-cube sampling tended to start at a bit less than double the optimized uncertainty. The bands roughly represent the spread of initial accuracies, and the fall of the bands to the optimal values represents the improvement one expects to gain through the optimization procedure.}
\end{minipage}
\end{figure}

The accuracy should increase with more training points. This is illustrated for the case of six model parameters with $\Lambda=3.0$ and $\alpha=0.01$ in Fig. \ref{fig:Sigma2vsNTrain}. The behavior underscores the value of choosing a number of training points corresponding to the number of independent coefficients for a linear, quadratic or cubic expansion. For six model parameters, those numbers of coefficients are $N_{\rm train}=7, ~28$ and 84 respectively. Once one of these thresholds is passed, the marginal value of an additional training point decreases. Then, as one approaches the next threshold, the marginal value increases again, though not to the degree of those points below the previous threshold.

The principal source of structure, or kinks, in Fig. \ref{fig:Sigma2vsNTrain} is due to whether the number of training points passes one of the thresholds mentioned above. However, there is a second less prominent source of structure that occurs when the number of training points corresponds to a symmetry. Given that this example was for Gaussian priors for each parameter, with the same widths (relative to $\Lambda$) for each prior, the prior had a spherical symmetry. Examples of arrangements with symmetries are simplexes and cubes. These might also include an extra point at the origin.  Such structures vary with dimensionality. For this case, $N_{\rm pars}=6$, the uncertainty for $N_{\rm train}=8$ appears noticeably less. On inspection, this turned out to be a simplex with an extra point at the center. Thus, one might well choose $N_{\rm train}=8$, even though it is one above the threshold. A smaller kink is seen at $N_{\rm train}=44$. But this structure is not so prominent as to drive one's decision to run more points.

For many cases, one finds that the minimization procedure might settle on a false minimum. The searches were initialized by either randomly placing the training points according to the prior weights, or using Latin hyper-cube samples \cite{lhs}. For Latin hyper-cube sampling each dimension of parameter space was divided into $N_{\rm train}$ equally probable slices. For each parameter, the value for each training point was set randomly forced to lie in a different slice. Figure \ref{fig:Sigma2vsNMC} shows how the accuracy metric improved with the number of random steps. Several chains are displayed to illustrate the progress of the search, and its dependence on the initial configuration. Settling on a false minimum can affect the estimated uncertainty by at most a few percent, and can lead to a bit of jitter for the optimized values. One can find the true minimum by repeating the minimization search several times, then choosing the best configuration. The number of random steps required to reach a desired accuracy increases with dimensionality of the model parameter space. The simple method applied here seems fine for cases where $N_{\rm pars}$ is near or less than a dozen, but if one had two dozen parameters, it might warrant a more sophisticated search algorithm.

Also shown in Figs. \ref{fig:Sigma2vsNTrain} and \ref{fig:Sigma2vsNMC} are the accuracy metrics for latin hyper-cube sampling. As this is not an optimized choice, there is more noticeable jitter in the curves. One can see that the optimization procedure reduces the expected uncertainty by nearly a factor of two at the quadratic threshold and by nearly a factor of five at the cubic threshold. Depending on the desired accuracy, one might need 50\% more runs with Latin hyper-cube sampling to achieve the same accuracy as would be obtained with an optimized placement. Figure \ref{fig:Sigma2vsNMC} also illustrates the improvement compared to random placement of the training points. In that case the training points were placed independently, according to the prior, with no correlation between points.

While performing the calculations of $\mathbf{E}_P(\sigma_Y^2)/\sigma_A^2$ considered here, the contributions to the uncertainty can be split into two pieces as expressed in Eq. (\ref{eq:sigma2prior}).  Whereas the first term in Eq. (\ref{eq:sigma2prior}) represented the uncertainty assuming $\Lambda$ is fixed at its most probable value, the second term represents the emulator uncertainty coming from $\Lambda$ varying, with the weight $Z$. Remarkably, this second contribution was never more than a few percent of the total expected uncertainty for the optimized placement, which makes it appear largely negligible. However, that fraction was often much larger when the optimization procedure was first started, perhaps a third of the uncertainty. Then as the training-point placement approached the optimal arrangement, this second contribution shrunk. It was observed that the first contribution, the part assuming $\Lambda$ was fixed, varied less with the training point placement compared to the reduction in the first contribution. Thus, the second contribution was rather instrumental in driving the placement of the training points, despite the fact that it is the lesser contribution. This would suggest that while this contribution is crucial for choosing the location of the training points, but that once that optimization is completed, and one has built an emulator, that contribution can be largely neglected when estimating the emulator's uncertainty. 

\subsection{Fine Structure of Point Placement}\label{subsec:finestructure}

Here, fine structure refers to the fact that certain numbers of training points might allow the placement to satisfy symmetries of the prior. If the prior is Gaussian in all directions, and is rotationally symmetric, then one might expect the optimum choice to exploit the symmetry. For the case of $N_{\rm train}=N_{\rm pars}+1$, one can place the training points in an $N_{\rm pars}$ dimensional simplex, a simplex being the placement of all points at the same radius, with the same distance between each pair. In two dimensions a simplex is an equilateral triangle, while for three dimensions the simplex is a tetrahedron. With this placement, there are a set of rotations under which each point can be rotated to another point in the structure without changing the structure. A few other arrangements may have similar symmetries. For example, one could have $N_{\rm pars}+1$ training points placed in a simplex, while one additional point might be placed in the center. For three dimensions, a buckyball configuration with $N_{\rm train}=60$ has a symmetry. But these symmetries are hard to come by, so in general it is not at all obvious how to place training points in an arbitrary number of dimensions, unless the number of training points happens to match some specific value. Further, even if the number of training points were to match that of a simplex, or of a simplex plus one extra point, one would have to determine the optimum radius of the simplex, which would necessitate a numerical search. Finally, in most cases, not all the model parameters have Gaussian priors, and if they were Gaussian, they might not have identical spreads. Another such symmetry is that where points are placed at the corners of an $N_{\rm pars}-$dimensional hyper-cube. This would require $2^{N_{\rm pars}}$ points. Again, another point might be added to the center. This arrangement would seem relevant for either the case where the prior was purely Gaussian with a spherical symmetry, or if all the priors were uniform with the same widths.

\begin{figure}
\centerline{\includegraphics[width=0.55\textwidth]{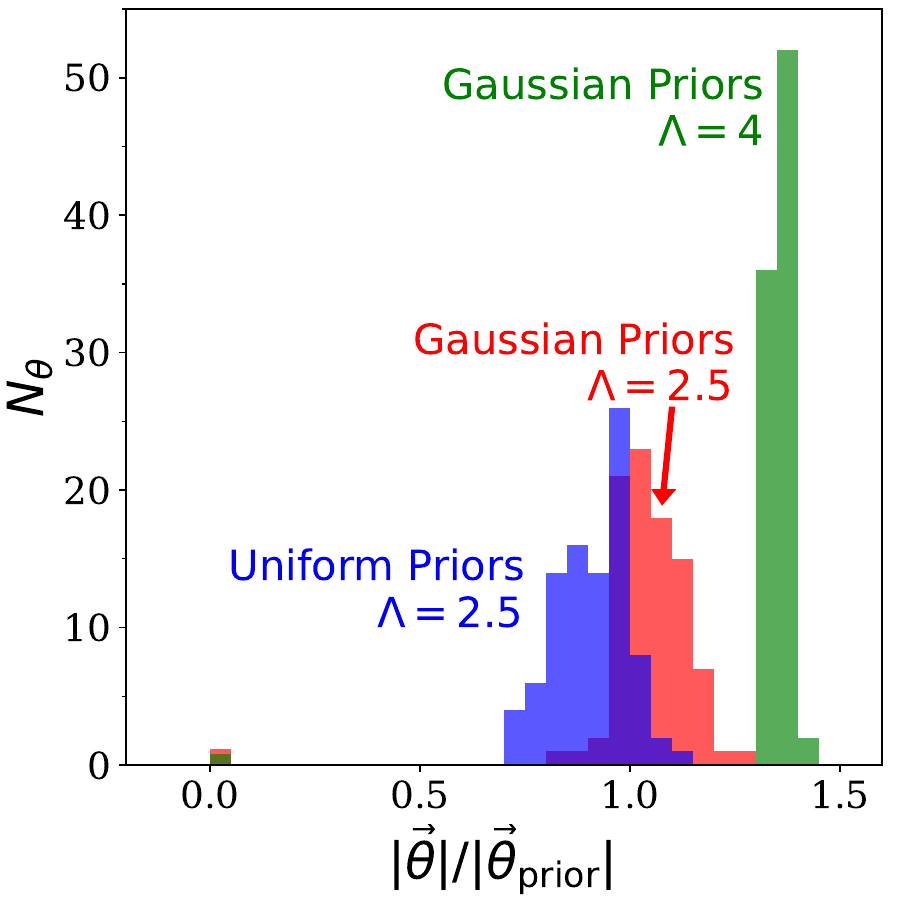}}
\caption{\label{fig:radii}
The distribution of the 91 optimized radii, $|\vec{\theta}_a|$, for the optimized placement of 91 training points in the 12-dimensional model-parameter space. For this case $N_{\rm train}=91$. If all model parameters had Gaussian priors with the same width relative to $\Lambda$, the radii tend to organize with one point in the center while the others have similar radius spread over different directions. For larger $\Lambda$ ($\Lambda=4$,green) larger radii were preferred compared to smaller $\Lambda$ ($\Lambda$=2.5,red). For uniform priors ($\Lambda=2.5$,blue) radii were smaller with a similarly broad distribution, and there was not training point in the middle.}
\end{figure} 
Thus, unless one has a symmetry, and additionally if the number of training points happens to allow placement according to that symmetry, the placement of the points is difficult to predict. The radius of a training point is simply the point's distance from the origin,
\begin{eqnarray}
|\vec{\theta}_a|&=&\left(\sum_{i<N_{\rm pars}}\theta_{a,i}^2\right)^{1/2}.
\end{eqnarray}
For the case with $N_{\rm pars}=12$ parameters, Fig. \ref{fig:radii} shows the radii of the various positions. For two cases, all the priors were Gaussian, with an r.m.s. radius of $1/\sqrt{3}$. The number of training points was set to $N_{\rm train}=91$, which corresponds to the quadratic-fit threshold for $N_{\rm pars}=12$. As can be seen the optimized placement is for most of the points to be a nearly the same radius, whereas one point is placed near the origin. When this same example was worked out for $N_{\rm pars}=6$, one point was at the center, while the others were all positioned at exactly the same radius, but for $N_{\rm pars}=12$ there was some small variation in the radii of the points away from the center. Two values of $\Lambda$ were considered in Fig. \ref{fig:radii}, $\Lambda=2.5$ and $\Lambda=4.0$. One can see that the radii moved outward for larger $\Lambda$.
Figure \ref{fig:radii} also shows the radii for 12 model parameters where each parameter had a uniform prior with minimum and maximum values of $-1$ and $+1$. That distribution is at somewhat smaller radius, and without a point placed at the origin.

\subsection{Overcoming the Pernicious Nature of High Dimensions}
\label{subsec:pernicious}

The overarching goal of this manuscript is to understand how to best emulate functions that can be described well with polynomials not extending to high order, or equivalently, are well described by Taylor expansions. In Sec. \ref{sec:theory} it was described how this can equivalently be posed in terms of Gaussian process approaches by requiring that the correlation length is larger than the size of the prior. Either way, emulation depended on defining a characteristic scale, $\Lambda$. For smaller $\Lambda$, more training points are needed. In high dimensionality  parameter spaces, it becomes more difficult to effectively cover the prior. In Fig. \ref{fig:sigma2vsNpars}, it was shown that even if the number of training points increased with number of quadratic coefficients, the accuracy of the emulator still fell as a function of $N_{\rm pars}$. The discussion here aims to provide several perspectives from which can appreciate the challenge of building high-dimensional emulators.

As discussed earlier and shown with Eq.s (\ref{eq:NcCounta}) and (\ref{eq:NcCountb}), the number of polynomial terms of order $n+1$ exceeds the number of terms of order $n$ by a factor proportional to the dimensionality of the parameter space, $N_{\rm pars}$. The variance due to the terms of order $n+1$ compared to those of order $n$ roughly grows a $N_{\rm pars}/\Lambda^2$. One can roughly state that the value of $\Lambda$ that is necessary to give a certain accuracy should increase as $\sqrt{N_{\rm pars}}$. Thus, accurate emulation will be more difficult in higher dimensions as illustrated in Fig. \ref{fig:sigma2vsNpars}. Another perspective from which one can gain the same conclusion is to consider the characteristic size of the parameter space. The characteristic r.m.s. sizes increase as $N_{\rm pars}^{1/2}$. The unconstrained covariance behaves as $e^{-(|\vec{\theta}_i-\vec{\theta}_j|)^2/2\Lambda^2}$, so to maintain the same characteristic correlations, $\Lambda^2$ must be increased proportional to $N_{\rm pars}$, or again, the $\Lambda$ must increase as $\sqrt{N_{\rm pars}}$. 

In higher dimensions it is increasingly critical that the convergence parameter be longer than the system. Otherwise, one would need to sample the function in each corner of the model-parameter phase space. For uniform priors, there are $2^{N_{\rm pars}}$ corners to the phase space. For 12 dimensions this requires 4096 training points, and for 18 dimensions, 262,144 training points. Thus, it can be untenable to place training points in each corner of the space. 

One can also see the problem with placing the training points in a spherically symmetric fashion. The hyper-volume of the parameter-space hyper-cube is $2^N$, whereas the volume of an $N-$dimensional hyper-sphere of radius $R=1$ is 
\begin{eqnarray}
V_{\rm sphere}&=&\Omega_N\int_0^R dr~r^{N-1}=\Omega_N\frac{R^N}{N}.
\end{eqnarray}
The solid angle, $\Omega_N$ in $N$ dimensions is
\begin{eqnarray}
\Omega_N&=&\frac{2\pi^{N/2}}{\Gamma(N/2)},
\end{eqnarray}
and after putting this together, the fraction of the hyper-cube's volume that is within the hyper-sphere is
\begin{eqnarray}
\frac{V_{\rm sphere}}{V_{\rm cube}}&=&\left\{
\begin{array}{rl}
\frac{(\pi/2)^{N/2}}{N!!},&N=2,4,6,8\cdots\\
\frac{(\pi/2)^{(N-1)/2}}{N!!},&N=3,5,7,\cdots\\
\end{array}
\right.
\end{eqnarray}
In two dimensions, the ratio is $\pi/4$, and in three dimensions it is $\pi/6$. In 10 dimensions it is $2.5\times 10^{-3}$. For high dimensions only a small fraction of the parameter space can ever lie inside a sphere used to place points. And, if the model is expensive, it may not be tenable to run the full model inside every corner. 

All the perspectives above emphasize the same conclusion. The only way a high-dimensional space can be well emulated is if the behavior is highly convergent. The parameter $\Lambda$ must be large compared to the characteristic size, which is $\sqrt{N_{\rm pars}/3}$. For 12 parameters this suggests $\Lambda>2$, and for 18 parameters this suggests $\Lambda>2.45$. If $\Lambda$ is near or below this limit, emulation may require a  large number of training points. 

For large numbers of model parameters it is crucial that one pause to consider whether some model parameters are more or less important. As stated in the beginning of Sec. \ref{sec:theory} this can be accomplished by scaling the model parameters. The model parameters with uniform priors can be scaled such that their priors have defined by $\theta_{{\rm max},i}=s_i$, where $s_i<1$ can be thought of as a suppression factor. The r.m.s. of the prior is then $s_i/\sqrt{3}$. The parameters with Gaussian priors can be scaled so the r.m.s. of their variance is $s_i/\sqrt{3}$. Scaling the range of the priors by $s_i$ is mathematically equivalent to scaling the convergence parameters, i.e. setting $\Lambda_i=\Lambda/s_i$. Scaling the priors, rather than scaling the convergence parameters, allows emulation to be simpler, as there is then only one convergence parameter. One could scale the parameters so that the most important parameters had $\theta_{{\rm max},i}=1$, or $\beta_i=3$. One would set $\Lambda$ to a value appropriate for those parameters, then scale $\theta_{{\rm max},i}$ and $\beta_i$ for the remaining parameters to account for the suppression of the expected sensitivity of $y(\vec{\theta})$ to the other parameters. 

If warranted, assigning suppression factors could be an effective means of dealing with high dimensionality. For example if $\Lambda=2$, it would be difficult to accurately emulate a 12-dimensional parameter space. However, if all but a couple of parameters were suppressed with $s_i=1/2$, the less critical parameters would not contribute much beyond their linear contribution. For example, if only three parameters were given full strength, then training-point placement might optimize so that the $2^3=8$ corners of the three-dimensional space were covered. This is much easier than attempting to cover $2^{12}=4096$ corners of parameter space. The software developed for this study was specifically designed to allow the user to assign suppression factors to each model parameter.

One might ask whether certain parameters might be deemed less important after the training data is created. That is certainly possible, and the emulator's uncertainty would fall whenever one would reduce $s_i$. However, assigning $s_i\ne 1$ also affects the optimized choice of training-point location. Thus, if the reduced importance is understood before the training data is generated that understanding can be exploited in the optimization of the training-point placement.

\setcounter{section}{3}

\section{Testing the Emulator and Training Point Optimization}\label{sec:testing}

\subsection{Overview of Testing}

For this study software was developed according to the methods described in the previous sections. The first piece of software generated a list of optimized training points as described in Sec. \ref{sec:optimize}. The input for this code is a list of model parameters and their priors, plus an estimate of the convergence parameter, $\Lambda$, and the point-by-point noise parameter $\alpha$. Each model parameter was assumed to have an independent prior, either a Gaussian form with a given centroid and width, or a uniform distribution between some minimum and maximum value. These parameters were then translated and scaled so that their priors were centered around zero.  By default, the widths were set according to the definitions in Eqs. (\ref{eq:uniformpriorsdef}) and (\ref{eq:gausspriorsdef}) so that the priors would have the same variance as a uniform distribution with $-1<\theta<1$ even if the prior distribution is Gaussian. An option was added to allow the width for a given parameter, $\theta_i$, to be scaled down by a factor $s_i$ if one believed a certain direction in parameter space was expected to contribute relatively less to the overall variance of the function. The software then provides a list of training points, $\vec{\theta}_a$, with $a=1-N_{\rm train}$. Some tests of this procedure were provided in the previous section. 

The second piece of software was the emulator discussed in Sec. \ref{sec:theory}. The emulator assumed the forms for $Y(\vec{\theta})$ and its uncertainty as described in Eq. (\ref{eq:summary}). This includes the contribution to the emulator uncertainty from the uncertainty in knowing $\Lambda$. The goal of this section is to study and test the emulator, and the theoretical foundations of both the emulator and the training-point optimization. The emulator's mean prediction, $Y(\vec{\theta})$, is exactly that of a typical GP emulator with a Gaussian kernel, despite the fact that the form was generated from a different perspective than the usual presentation of GP emulators. Thus, the tests here will focus on the following questions.
\begin{enumerate}\itemsep=0pt
\item If a model is built based on a polynomial expansion with specific values of $\Lambda$ and $\sigma_A$, how well will the tuned emulator reproduce those values?
\item In the limits of $\Lambda\gtrsim 2$, can one not only reproduce a model, but can it also provide a meaningful estimate of its uncertainty?
\item Is the contribution to the emulator uncertainty from the uncertainty in $\Lambda$ significant?
\item Did the estimate of the emulator's uncertainty averaged over the prior, $\mathbf{E}_P(\sigma_Y^2)/\sigma_A^2$, from the training point selection before any model runs were processed, match that of the tuned emulator?
\end{enumerate}
A subsection is devoted to each of the four questions above.

In order to test the models, two model types were created. The first type, Type-$A$, is based on creating polynomials with random coefficients, where the variance of the coefficients is consistent with the assumed form of the emulator for assumed values of $\sigma_A$ and $\Lambda$. The second, Type-$B$, is a sum of randomized trigonometric and exponential functions. Test functions were thus generated of both types, and for each function an emulator was constructed. By creating thousands of such models, and comparing emulators trained upon them, one can test how accurately the emulator represents its uncertainty. One can also test whether the estimates of accuracy used for the optimization of training-point placement are consistent with the accuracy actually exhibited by the trained emulator. 

\subsubsection{Type-$A$ Model}
The form of the expansion is the same as what was assumed when designing the emulator:
\begin{eqnarray}
F(\vec{\theta})&=&e^{-|\vec{\theta}|^2/2\Lambda^2}\sum_{\vec{n}}\frac{A_{\vec{n}}}{\sqrt{n_1!n_2!\cdots n_{N_{\rm pars}}}!}
\left(\frac{\theta_1}{\Lambda}\right)^{n_1}\left(\frac{\theta_2}{\Lambda}\right)^{n_2}\cdots\left(\frac{\theta_{N_{\rm pars}}}{\Lambda}\right)^{n_{N_{\rm pars}}}.
\end{eqnarray}
Here, the coefficients $A_{\vec{n}}$ are assumed to have a Gaussian distribution with variance $\sigma_A^2$. To create the model, one chooses the dimensionality, $N_{\rm pars}$, the variance $\sigma^2_A$ and the convergence parameter $\Lambda$. The coefficients $A_{\vec{n}}$ are then randomly set according to a Gaussian distribution, with the maximum order set to $K_{\rm max}=5$. Because this is the same form assumed by the emulator, aside from limiting $K\le K_{\rm max}$, one can test  whether the emulator, $Y(\vec{\theta})$, will choose optimum values of the hyper-parameters, $\sigma_A$ and $\Lambda$, that match those used to create the full model, $F(\vec{\theta})$. A different choice for the coefficients represented a different model. For this study thousands of Type-$A$ models were considered.

\subsubsection{Type-$B$ Model}
The second functional form was constructed by summing over randomly generated sums of trigonometric and hyperbolic functions. The observable was a sum of six functions, each constructed from six model parameters.
\begin{eqnarray}\label{eq:randomeqs}
F(\vec{\theta})&=&D_0\cos(t_0/\Lambda_0)+D_1\sin(t_1/\Lambda_1)+D_2\tan^{-1}(t_2/\Lambda_2)
+D_3\cosh(t_3/\Lambda_3)+D_4\sinh(t_3/\Lambda_4)+D_5\tanh^{-1}(s_5),\\
\nonumber
s_5&\equiv&\frac{t_5}{\sqrt{\Lambda_5^2+t_5^2}},\\
t_i&=&\sum_j r_{ij}\theta_j/\sqrt{6}+r_{ijk}\theta_j\theta_k/\lambda_{jk}.
\end{eqnarray}
Here, $r_i$ and $r_{ij}$ and $r_{ijk}$ are random numbers taken from a uniform distribution between $\pm 1$. The values $\Lambda_i$ and $\lambda_{jk}$ were random numbers between 2.0 and 5.0. Each of the coefficients $D_i$ were randomly chosen according to a Gaussian distribution with the same variance. Again, thousands of such models were generated to perform tests.

\subsection{Choosing the Hyper-Parameters}

In order to test the algorithms from Sec. (\ref{sec:theory}) for choosing the hyper parameters, $\sigma_A$ and $\Lambda$, $10^5$ Type-$A$ models were constructed. Each model was constructed using $\sigma_A=100$ and $\Lambda=3$. Because the type-$A$ models use the form assumed by the emulator, one might expect the algorithm to choose values of $\sigma_A$ and $\Lambda$ that match the values used to generate the models. There should be some fluctuation due to the coefficients for the models having been chosen with random distributions. 
\begin{figure}[ht]
\centerline{\includegraphics[width=0.6\textwidth]{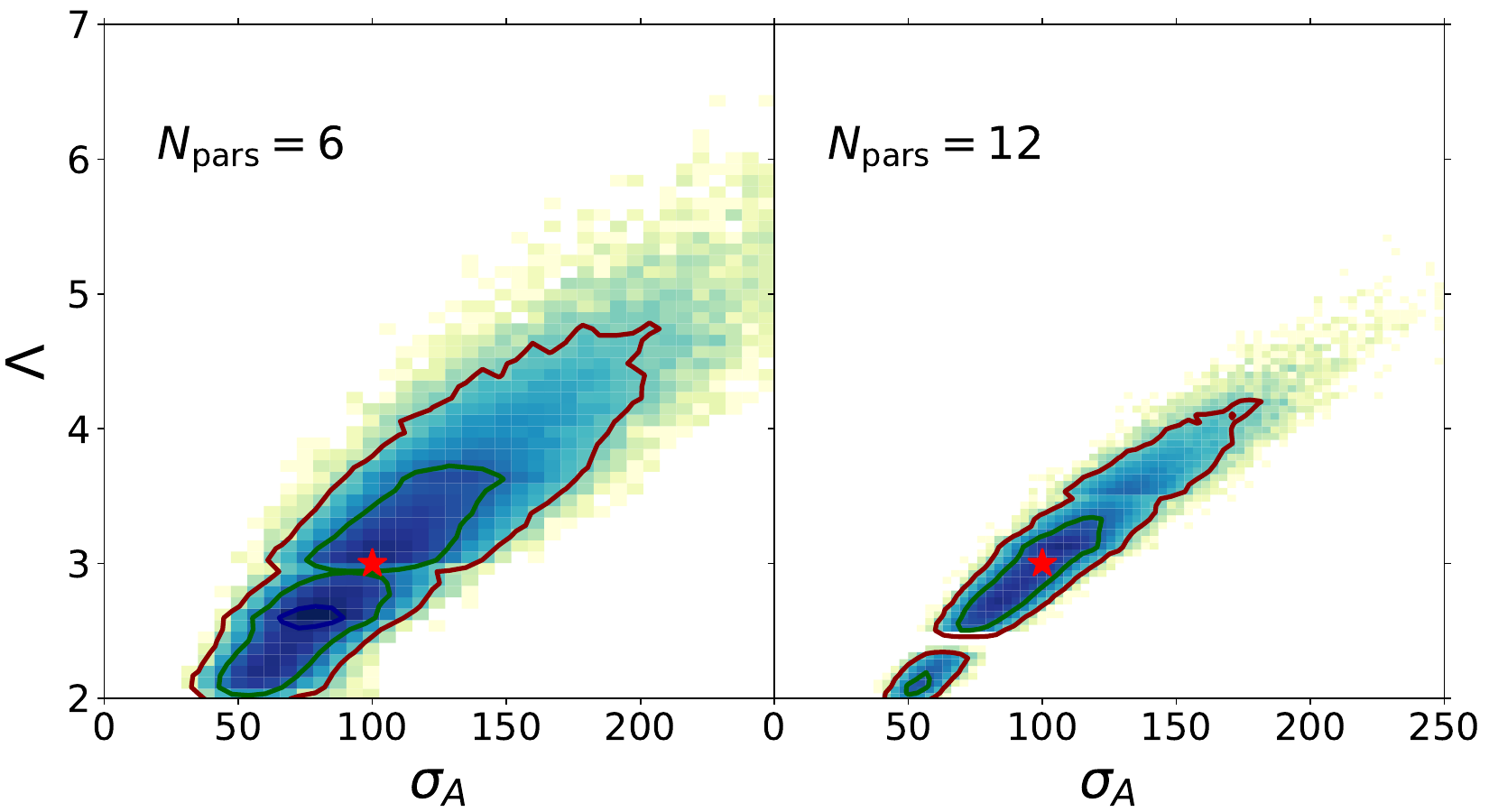}}
\caption{\label{fig:sigmalambda}
$10^5$ models were constructed based on a expansion with random coefficients chosen according to the form assumed by the emulator using $\Lambda=3$ and $\sigma_A=100$. The hyper parameters chosen by the tuning methods from Sec. \ref{sec:choosehyper} are compared to the values used to build the model (red stars). For each model the extracted hyper-parameters varied significantly, but if one averages the extracted hyper-parameters over the $10^5$ models, the average values come within a percent of $\Lambda=3$ and $\sigma_A=100$. A strong band structure was found in that if the extracted value of $\sigma_A$ was low, the value of $\Lambda$ also tended to be low. Both panels exhibit some  bimodal behavior in the distribution.}
\end{figure}
The distribution of extracted hyper parameters for the $10^5$ models is displayed in Fig. \ref{fig:sigmalambda}. Distributions were found for two cases, one with 6-parameter models, and one with 12-parameter models. For the 6-parameter models, the emulators were trained at 28 points in parameter space, and for the 12-parameter model tuning involved 91 training points. These numbers correspond to the number coefficients that have $K_{\rm max}\le2$. On average, the hyper parameters well matched those of the model (shown by the red stars), $\sigma_A=100$ and $\Lambda=3$. Averaging over the distribution, the extracted hyper parameters matched those used to generate the model to within 1\%. However, the distribution was rather wide, with widths of a few tens of percent.  A strong correlation was observed between the extracted values of $\sigma_A$ and $\Lambda$, as one can see a strong band structure. Larger $\Lambda$ is highly correlated with larger $\sigma_A$. For both cases, the distributions exhibit some bimodal behavior. Being that the emulator uncertainty is proportional to $\sigma_A$ and falls with increasing $\Lambda$, the overall uncertainty did not vary greatly from one instance to another despite the rather broad range of the distribution. The variance of the extracted hyper parameters decreased for the case with larger numbers of model parameters.

If one extracts small values of $\Lambda$, the expansion is poorly converging. For that reason, whenever the extracted value fell below $\Lambda=\sqrt{N_{\rm pars}/3}$, the emulator would set $\Lambda$ to that value. This was a rare occurrence for this case.  A second issue can arise for large $\Lambda$, especially for larger numbers of parameters. The tuning procedure involves inverting the matrix $B$ and finding its determinant. The form for $B$ is
\begin{eqnarray}
B_{ab}&=&C_0(\vec{\theta}_a,\vec{\theta}_b) +\Delta_{ab},\\
\nonumber
\Delta_{ab}&=&\alpha^2\delta_{ab}.
\end{eqnarray}
The term with $\vec{n}=0$ provides a contribution of unity to each element of the $N_{\rm pars}\times N_{\rm pars}$ matrix. Thus, ignoring the other contributions, one would have $B_{ab}=1$, and $B$ would be highly singular. Contributions from other $\vec{n}$ behave as $\Lambda^{-|\vec{n}|}$, and keep the matrix from being singular, but numerically, overcoming the nearness to the singularity is more difficult for large $\Lambda$ or for higher dimensionality of the $N_{\rm pars}\times N_{\rm pars}$ matrix. Further, some choices of training points can accidentally lead to singular, or nearly singular, matrices. Fortunately, adding point-by-point noise, $\alpha\ne 0$, also eliminates the singularity. The same issues can plague GP emulators. In that case the noise term above is referred to as a ``nugget'' \cite{RasmussenGP}. 

\subsection{Predicting The Uncertainty}

\begin{figure}[ht]
\centerline{\includegraphics[width=0.75\textwidth]{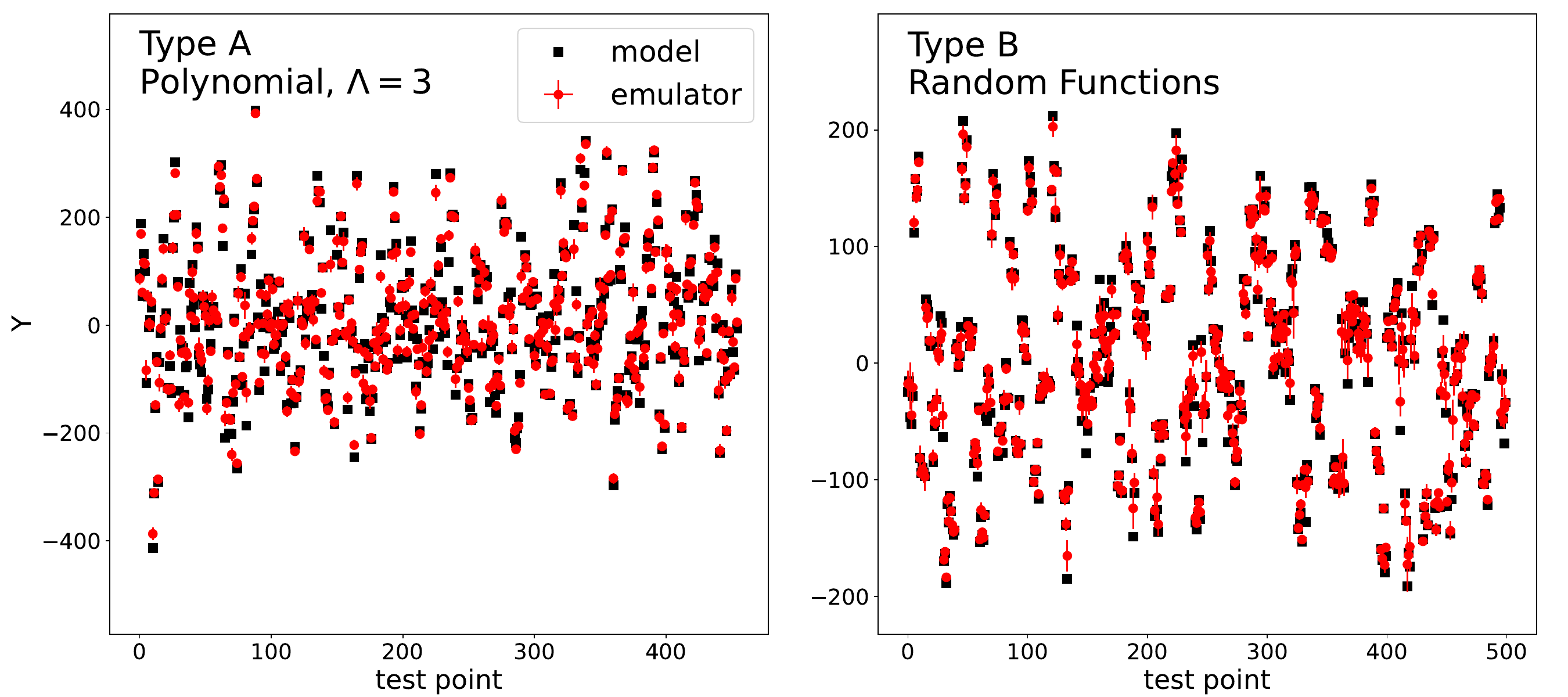}}
\caption{\label{fig:YvsY}
Emulator predictions and their uncertainties (red circles) are compared to model calculations (black squares) at 500 random points in parameter space.  Models involved 12 parameters and were trained at 91 points. For the left-side panel, the model was a polynomial with random coefficients, type-$A$, chosen with the emulator's assumed form with $\Lambda=3$, while in the right-side panel the model is a sum of random assorted trigonometric and exponential functions, type-$B$ as described by Eq. (\ref{eq:randomeqs}). }
\end{figure}

To address questions (2) and (3) from the beginning of this section, both Type-$A$ and Type-$B$ models were applied. For each model emulators were trained from model values at the training points. The training points were chosen according to method outlined in Sec. \ref{sec:optimize}. For this study, the  parameter space was 12-dimensional with Gaussian priors. The training involved evaluating the model at 91 optimized points, 91 being the threshold for a quadratic fit.  Figure \ref{fig:YvsY} presents this comparison for 500 of these points from a variety of different models, and for a variety of random points in model-parameter space. The emulator's uncertainty for each point is also displayed. The percentage of emulations within one unit of its stated uncertainty should be approximately 68.3\% for Gaussian errors. For the Type-$A$ models the success rate was indeed 68\%, while for the random assorted functions, the success rate was 88\%.  
\begin{figure}
\centerline{\includegraphics[width=0.35\textwidth]{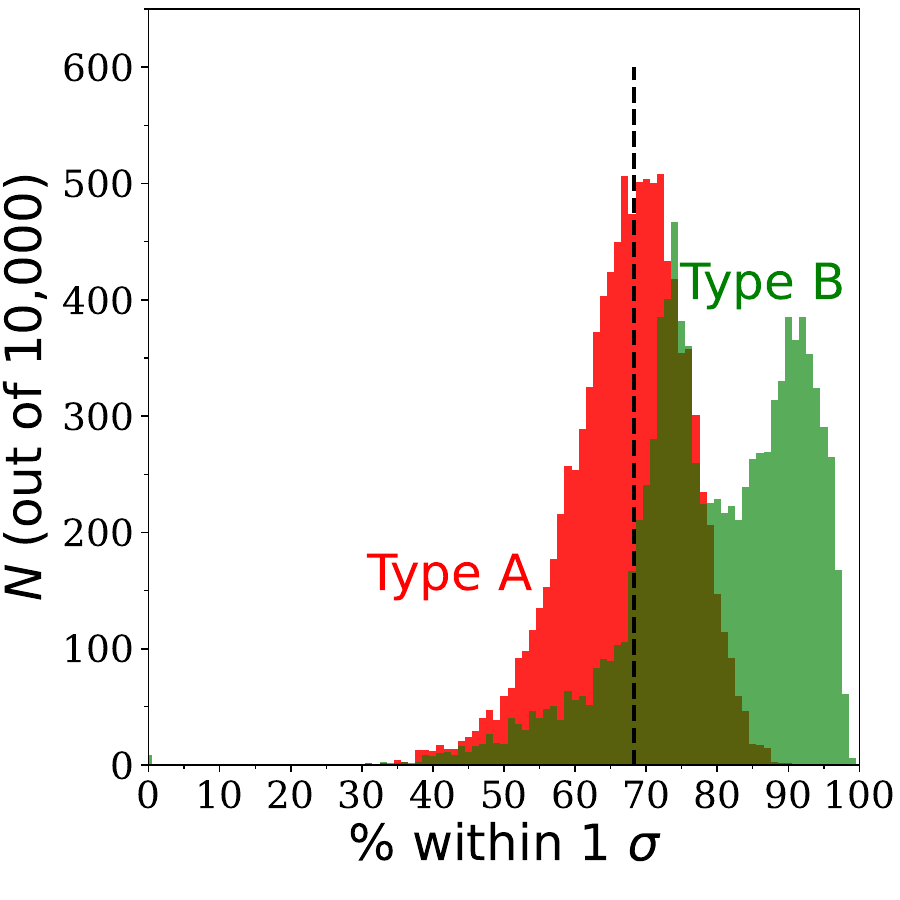}}
\caption{\label{fig:fitpercent}
Ten thousand models were created, each a polynomial expansion with random coefficients chosen to match the form assumed by the emulator. After the emulator was tuned, each model was evaluated at 1000 random points within the prior and compared to the emulator prediction. For each model the percentage was found for the prediction being within $\pm \sigma_Y$. If the emulator were to perfectly estimate its uncertainty, the distribution would be peaked around 68.3\%. For the models of type $A$ the emulator tended to predict its uncertainty rather well, whereas for type $B$ the emulator tended to modestly overpredict its uncertainty.}
\end{figure}

When calculating the uncertainties at the 500 random points used to generate Fig. \ref{fig:YvsY}, the contribution to the uncertainty from the uncertainty in determining $\Lambda$ was included. That contribution, related to term proportional to $W_{\Lambda\Lambda}$ in Eq. (\ref{eq:summary}) was rarely more than a few percent of the total uncertainty. This was consistent with the expectations for that contribution, which were  calculated when the choice of training points was being optimized. 

\subsection{Testing the Preconstrained Estimates of the Uncertainty}
\begin{figure}
\centerline{\includegraphics[width=0.45\textwidth]{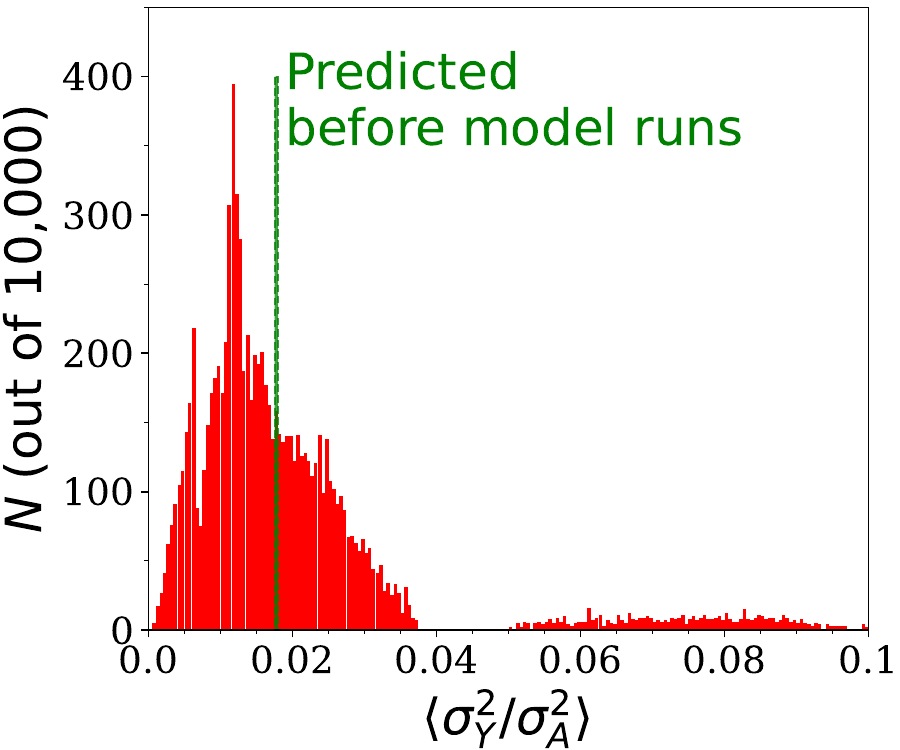}}
\caption{\label{fig:sigma2count}
When optimizing the choice of training points, the procedure produced estimates of the prior-averaged uncertainty relative to $\sigma_A^2$. That value is represented by the green dashed line. From tuning 10,000 randomly generated type-$A$ models with $N_{\rm pars}=12$, a distribution of the accuracy metrics, $\mathbf{E}_P(\sigma_Y^2)/\sigma_A^2$, is displayed here. The distribution is remarkably wide, but the average indeed matches the description. The large width mainly comes from the fact that the value of $\Lambda$ varies, and that the estimated uncertainty is strongly sensitive to $\Lambda$. 
}
\end{figure}

While choosing the training points, and before any full-model runs were performed, estimates were made of the expected emulator's uncertainty averaged over the prior. Figure \ref{fig:sigma2count} illustrates the degree to which that prediction was successful. Ten thousand models of type-$A$ were randomly generated using $N_{\rm pars}=12$ with Gaussian priors. The coefficients were random but chosen with weights assuming $\Lambda=3$. The distribution's average indeed matched the predicted value, but the width of the distribution is remarkably wide. The bimodality seen in Fig. \ref{fig:sigmalambda} is also manifest here. An island of fits with smaller $\sigma_A$ and smaller $\Lambda$ in Fig. \ref{fig:sigmalambda} is seen in Fig. \ref{fig:sigma2count} as an additional island of values of larger $(\sigma_Y/\sigma_A)^2$. Even though these models tended to have larger values of $\sigma_Y/\sigma_A$, their overall uncertainties were similar to those from the bulk of the distribution because $\sigma_A$ tended to be smaller. Hence, the distribution of accurate predictions in Fig. \ref{fig:fitpercent} was not particularly wide despite the widths of the distributions in Fig. \ref{fig:sigmalambda} and \ref{fig:sigma2count}.

\setcounter{section}{4}
\section{Summary}

The goal of this study was to consider the emulation of functions which are expected to be well described by polynomials with a few orders in the expansion. From this perspective, questions were addressed about what forms of the emulator one should use, how one should choose the training points, and whether one can trust the emulator's statement of its uncertainty. In particular, the motivation is driven by the need to consider models with high-dimensional parameter spaces that, due to their numerical cost, should be trained with optimum efficiency. 

In Sec. \ref{sec:theory} it was found that if the following criteria were set, then the form of the emulator would be that of a Gaussian process with a Gaussian kernel. The criteria were:
\begin{enumerate}\itemsep=0pt
\item The function should be analytic throughout the prior.
\item At any point the function should fluctuate according to a Gaussian distribution.
\item The prior-to-constraint covariance should depend only on $|\vec{\theta}_1-\vec{\theta}_2|^2$.
\end{enumerate}
The polynomial form corresponding to these criteria was found to be:
\begin{eqnarray}\label{eq:summary_Aexpansion}
Y(\vec{\theta})&=&\sum_{\vec{n}}^{N_{\rm coeff}}A_{\vec{n}}\mathcal{T}_{\vec{n}}(\vec{\theta}),\\
\nonumber
\mathcal{T}_{\vec{n}}(\vec{\theta})&=&
\frac{1}{\sqrt{n_1!n_2!\cdots n_{N{\rm pars}}!}}
\left(\frac{\theta_1}{\Lambda}\right)^{n_1}\left(\frac{\theta_2}{\Lambda}\right)^{n_2}\cdots 
\left(\frac{\theta_N}{\Lambda}\right)^{n_{N{\rm pars}}}e^{-|\vec{\theta}|^2/2\Lambda^2}.
\end{eqnarray}
Expressions were found for all expansion coefficients in terms of the training data. The actual coefficients were not needed given that the emulator predictions and uncertainties could all be expressed in terms of the unconstrained covariance,
\begin{eqnarray}
\Cov_0(Y(\vec{\theta}_1),Y(\vec{\theta}_2))&=&\sigma_A^2e^{-|\vec{\theta}_1-\vec{\theta}_2|^2/2\Lambda^2}.
\end{eqnarray}
In practice, the emulator behaves exactly as a GP emulator with a Gaussian kernel, and with the convergence parameter $\Lambda$ being the same as the correlation-length hyper-parameter in the GP paradigm. Whereas a large convergence parameter, $\Lambda$, relative to the prior size, forces the polynomial to converge rapidly, within the GP paradigm one can equivalently state that the correlation length, $\Lambda$, should well exceed the size of the prior.

The emulator's uncertainty was shown to have two contributions. The first comes from the fact that the expansion coefficients, $A_{\vec{n}}$, are constrained by the training but are still very much unknown given that there are typically many more expansion coefficients than training constraints. This first contribution is identical to that which is typically used for GP emulators. The second contribution comes from the fact that the hyper-parameters, $\sigma_A$ and $\Lambda$, are not precisely known. Although there are expressions for the most likely values, the likelihood for varying $\sigma$ and $\Lambda$ is characterized by widths. Because the width $W_{\Lambda\Lambda}\ne 0$ and because the prediction can be sensitive to $\Lambda$, this leads to a second contribution to the emulator uncertainty. For cases with reasonably large numbers of training points, this second contribution turned out to be relatively small but played an outsized role in the optimization of training-point placement.

In Sec. \ref{sec:optimize}, a procedure was presented for optimizing the placement of training points. The metric for accuracy was the squared uncertainty averaged over the prior, $\mathbf{E}_P(\sigma_Y^2)/\sigma_A^2$. Aside from the location of the training points, the accuracy metric also depended on a guess for what value of $\Lambda$ might be found once the actual training was performed. The measure was expressed analytically for Gaussian or uniform priors. Even with a rather inefficient search procedure, optimized configurations of training points were readily found, although for very high dimensions the procedure might take hours. The optimization was driven mainly by the second source of uncertainty, that due to $\Lambda$ varying. However, once an optimized placement was chosen, this second source of uncertainty was found to be rather negligible.

It was found that the optimized configuration would typically improve on simple Latin hyper-cube sampling by a factor of nearly two or more as quantified by the accuracy metric. This approach might thus reduce the number of full-model training runs by perhaps a third. The marginal improvement of an additional full-model run was found to fall when the number of training points passed thresholds related to the number of expansion coefficients of a given order. The lesson gained here is that it could be advantageous to choose a number of full-model runs that exactly correspond to the number of coefficients up to a certain order. For example, for $N_{\rm pars}=12$, there are 91 coefficients of quadratic order or less. The expected accuracy improved significantly more while increasing the number of runs to 91, and improved less with additional runs. 

The measure of expected accuracy, $\mathbf{E}_P(\sigma_Y^2)/\sigma_A^2$, was found to be highly sensitive to the convergence parameter, $\Lambda$ and to the dimensionality of the model-parameter space. Sensitivity to the point-by-point noise was small as long as the noise was small compared to the other variation of the model. Although the estimate of the accuracy depended strongly on $\Lambda$, the actual optimum locations of the training points changed only modestly with $\Lambda$. The sensitivity to the dimensionality is a prime lesson to be taken from this study. Even if the number of training points is increased to match the increased number of coefficient of a given order as $N_{\rm pars}$ increases, the expected accuracy falls. The challenge of emulating high-dimension spaces was discussed in detail at the end of Sec. \ref{sec:optimize}.

To test the methods discussed above, two simple models were considered. The first was a model based on the expansion in Eq. (\ref{eq:summary_Aexpansion}), using random expansion coefficients. This precisely coincides with the form assumed by the emulator. The second form was a sum of random trigonometric and exponential forms. Thousands of models were constructed for each class. The models were each then used to train and build emulators. Emulator predictions were then compared to model values at randomly selected points throughout the prior. For the first class of models the emulator provided an outstanding estimate of its accuracy. This was expected given that the model's form matched that assumed by the emulator. For the case of randomized trigonometric and exponential functions, the emulator modestly understated its accuracy. 

For the first type of model, one would expect the emulator's value of $\sigma_A$ and $\Lambda$ to match those used to define the model. On average the extracted values matched well, but the distribution of extracted values was large. Despite the fact that $\sigma_A$ and $\Lambda$ might stray far from the values used to generate the functional data, and despite the fact that the ratio $\sigma_Y/\sigma_A$ often strayed far from the predictions provided by the training-point-placement procedure, the emulators largely provided good estimates of their uncertainty. This was driven by the fact that extracted $\sigma_A$ and $\Lambda$ values were highly correlated. Thus, if the extracted $\sigma_A$ was much higher than the value used to generate the training data, the extracted values of $\Lambda$ were also quite high. The two taken together, then led to rather modestly varying uncertainties. 

An important lesson to taken from this study is that high-dimension parameter spaces are difficult to emulate unless the functions are highly convergent. However, in most instances one expects many of the model parameters to contribute weakly. A strategy was provided in Sec. \ref{subsec:pernicious} to exploit this fact by simply assigning suppression factors to some subset of model parameters. The suppression factors  would allow the procedure to account for the relatively weaker sensitivity to some parameters by scaling their priors to narrower ranges, while assuming all parameters contributed according to a single convergence parameter. 

Aside from the lessons gained from this study, the work points to places where the procedures could be improved or expanded. The main opportunity for improvement would be in coding a more efficient procedure for optimizing the training points. Steepest descent algorithms might accelerate the search by multiple factors of two. A primary opportunity for expanding the functionality would be to consider more classes of priors when optimizing the training-point placement. The current software considered Gaussian and uniform priors. One could consider half Gaussians, or even or odd functions. 

Finally, one needs to ask the question as to whether one should perform full model runs simultaneously or sequentially. Here, by doing a significant number simultaneously, one can maximize the accuracy for that number of runs as measured over the entire prior. However, once the emulator is tuned well enough, one could perform a Markov Chain Mote Carlo (MCMC) exploration of space and identify those regions of the parameter space which are likely. One could then augment the training data with model runs specifically aimed to best constrain predictions for regions that are relevant to the goal of determining the model parameters. The additional points would likely be located in the likely regions of parameter space. Given that one has little concrete knowledge of where the likely region of parameter space lies before the training, it might make sense to simultaneously perform a first set of predictions to constrain the linear or quadratic behavior without comparing to experiment. Then, after determining the likely regions of parameter space by comparing to data, one might perform additional runs in a more sequential manner. Unfortunately, it is difficulty to say whether the first set of runs should be small, e.g. only enough for a linear fit, or whether one might opt for enough for a quadratic fit. The decision might rest on numerous factors. This might include how convergent the functions were expected to be, their numerical cost, and whether one can easily perform some number of full-model runs simultaneously. Unfortunately, their convergence would not be well understood until after the first round of model runs. Choosing an optimized strategy for most cases should be quite beneficial and allow one to avoid performing unnecessary full-model runs.

\appendix
\section{Properties of the Polynomial Basis Functions}\label{sec:basisprops}

The polynomial expansion presented in Sec. \ref{sec:theory} can be considered as an expansion in basis functions. The purpose of this appendix is to examine the orthonormality and completeness properties of those basis functions.

The polynomial expansion, for $N=N_{\rm pars}$, used for the studies here is of the form
\begin{eqnarray}
Y(\vec{\theta})&=&e^{-|\vec{\theta}|^2/2\Lambda^2}\sum_{\vec{n}}\frac{A_{\vec{n}}}{\sqrt{n_1!n_2!\cdots n_{N}}!}
\left(\frac{\theta_1}{\Lambda}\right)^{n_1}\left(\frac{\theta_2}{\Lambda}\right)^{n_2}\cdots\left(\frac{\theta_{N}}{\Lambda}\right)^{n_{N}},
\end{eqnarray}
where $A_{\vec{n}}$ were independent coefficients, each obeying a Gaussian distribution,
\begin{eqnarray}
W(A_{\vec{m}})&=&\frac{1}{\sqrt{2\pi\sigma_A^2}}e^{-A_{\vec{m}}^2/2\sigma_A^2},\\
\nonumber
E_0 (A_{\vec{m}},A_{\vec{n}})&=&\sigma_A^2\delta_{\vec{m}\vec{n}}.
\end{eqnarray}
This can be considered as an expansion in terms of basis functions, $\phi_{\vec{n}}$ defined as
\begin{eqnarray}
\phi_{\vec{n}}(\vec{\theta})&\equiv&\frac{1}{(2\pi\Lambda^2)^{N/4}}
e^{-|\vec{\theta}|^2/2\Lambda^2}
\frac{1}{\sqrt{n_1!n_2!\cdots n_{N}}!}
\left(\frac{\theta_1}{\Lambda}\right)^{n_1}\left(\frac{\theta_2}{\Lambda}\right)^{n_2}\cdots\left(\frac{\theta_{N}}{\Lambda}\right)^{n_{N}}.
\end{eqnarray}
Here, the basis functions are labeled by $\vec{n}$, which is a list of any $N$ non-negative integers. With this definition the emulator is
\begin{eqnarray}
Y(\vec{\theta})&=&(2\pi\Lambda^2)^{N/4}\sum_{\vec{n}}A_{\vec{n}}\phi_{\vec{n}}(\vec{\theta}).
\end{eqnarray}
The prefactor, $1/(2\pi\Lambda^2)^{N/4}$, is defined for convenience.

First, the orthogonality and normalization properties are:
\begin{eqnarray}
Z_{\vec{m}\vec{n}}&=&\int d\theta_1\cdots d\theta_{N}~\phi_{\vec{m}}(\vec{\theta})\phi_{\vec{n}}(\vec{\theta})\\
\nonumber
&=&\prod_{j=1}^{N}z_{m_jn_j},\\
\nonumber
z_{m_jn_j}&=&\frac{1}{(2\pi\Lambda^2)^{1/2}}\frac{1}{\sqrt{m_j!n_j!}}
\int d\theta~e^{-\theta^2/\Lambda^2}
\left(\frac{\theta}{\Lambda}\right)^{m_j+n_j}\\
\nonumber
&=&\left\{\begin{array}{rl}
\frac{1}{(2\pi)^{1/2}}\frac{1}{\sqrt{m_j!n_j!}}\Gamma([m_j+n_j+1]/2),&m_j+n_j~{\rm is~even}\\
0,&m_j+n_j~{\rm is~odd}
\end{array}\right. ~.
\end{eqnarray}
The basis functions are neither normalized nor orthogonal, although the normalizations are of order unity.  For the purposes of emulation it is not clear that orthonormality properties play an important role. 

On the other hand, understanding completeness determines whether the basis can capture all facets of the function it is attempting to represent. The completeness properties are defined by
\begin{eqnarray}
\sum_{\vec{n}}~\phi_{\vec{n}}(\vec{\theta}_a)\phi_{\vec{n}}(\vec{\theta}_b)&=&
\frac{1}{(2\pi\Lambda^2)^{N/2}}e^{-|\vec{\theta}_a|^2/2\Lambda^2}e^{-|\vec{\theta}_b|^2/2\Lambda^2}\\
\nonumber
&&\sum_{\vec{n}}\frac{1}{n_1!\cdots n_{N!}}
\left(\frac{(\theta_{a1}\theta_{b1})^{n_1}}{\Lambda^{2n_1}}\right)\cdots
\left(\frac{(\theta_{aN}\theta_{bN})^{n_N}}{\Lambda^{2n_N}}\right)\\
\nonumber
&=&\frac{1}{(2\pi\Lambda^2)^{N/2}}e^{-|\vec{\theta}_a-\vec{\theta}_b|^2/2\Lambda^2}.
\end{eqnarray}
For a complete basis the right-hand side would have been a delta function. Instead, it is a smeared delta function. It integrates to unity but has a non-zero width, $\Lambda$. This would satisfy the completeness relation if $\Lambda$ were allowed to approach zero, but in that limit one would likely require an infinite number of terms in the expansion. 

This basis can then be considered as being incomplete in that it cannot well represent small scale features, i.e. those with characteristic scales below $\Lambda$. For perfect completeness to hold over an infinite volume one would expect to require a continuous set of coefficients that cover an infinite extent. A similar example would be a Fourier transform, where instead of the infinite sets of integers, $\vec{n}$, one has the vector $\vec{k}$ which is distributed continuously. If one confined the vectors $\vec{k}$ to integer multiples of some step vectors, $\Delta k_i$, one would be unable to represent features below the inverse size of the step vectors, and if one cut off the expansion at some maximum, $k_{\rm max}$, one would be unable to describe structures of sizes below $1/k_{\rm max}$. For the problem here, $\Lambda$ sets the minimum scale at which one would expect features for the full function, and the maximum scale at which one can represent the function is determined by whether the sum over $\vec{n}$ is cut off at some maximum order.

Above, for a given rank $K$, the basis functions include independent terms for all $n_1+n_2+\cdots n_{N_{\rm pars}}=K$. One might have instead chosen to express the polynomials as all terms $i_1,i_2,i_K$, so that $Y(\vec{\theta})$ would then be:
\begin{eqnarray}\label{eq:alternateform}
    Y(\vec{\theta})&=&\sum_{K,i_1\cdots i_K}A_{i_1\cdots i_K}\frac{1}{K!}
    \theta_{i_1}\theta_{i_d}\cdots\theta_{i_K}.
\end{eqnarray}
Here, $1\le i_j\le N_{\rm pars}$. In this form, there are more coefficients. In the first form there was a term proportional to $\theta_1^2\theta_2$, and in this form there are three separate terms with the same product, $\theta_1\theta_1\theta_2$, $\theta_2\theta_1\theta_1$ and $\theta_1\theta_2\theta_1$. This alternate form has redundant coefficients, but those coefficients have less strength. In the previous form $N!d_{\vec{n}}^2$ was set to $N!/(n_1!n_2!n_3...)$, which equals the number of identical terms contributing in the new form. Because having $d$ independent gaussians is the same as having a single gaussian of strength $\sqrt{d}$, the end result is the same correlation. Thus, the alternate form in Eq. (\ref{eq:alternateform}) leads to all the same results as that presented in Sec. (\ref{sec:theory}).

\begin{acknowledgments}
This work was supported by the Department of Energy Office of Science through grant number DE-FG02-03ER41259 and by the National Science Foundation's  Cyberinfrastructure for Sustained Scientific Innovation program through grant number OAC-2004601. The authors express gratitude to Simon Mak and Daniel Phillips for generously sharing their insight into the subject.
\end{acknowledgments}

\setcounter{section}{5}


\begin{thebibliography}{99}

\bibitem{Domingues:2024pom}
T.~S.~Domingues, R.~Krupczak, J.~Noronha, T.~N.~da Silva, J.~F.~Paquet and M.~Luzum,
Phys. Rev. C \textbf{110}, no.6, 064904 (2024).

\bibitem{JETSCAPE:2024cqe}
R.~Ehlers \textit{et al.} [JETSCAPE],
[arXiv:2408.08247 [hep-ph]].

\bibitem{JETSCAPE:2024ofm}
R.~Ehlers \textit{et al.} [JETSCAPE],
EPJ Web Conf. \textbf{296}, 15009 (2024).

\bibitem{Paquet:2023rfd}
J.~F.~Paquet,
J. Phys. G \textbf{51}, no.10, 103001 (2024).

\bibitem{JETSCAPE:2023ikg}
W.~Fan \textit{et al.} [JETSCAPE],
Phys. Rev. C \textbf{109}, no.6, 064903 (2024).

\bibitem{Heffernan:2023gye}
M.~R.~Heffernan, C.~Gale, S.~Jeon and J.~F.~Paquet,
Phys. Rev. Lett. \textbf{132}, no.25, 252301 (2024).

\bibitem{Heffernan:2023utr}
M.~R.~Heffernan, C.~Gale, S.~Jeon and J.~F.~Paquet,
Phys. Rev. C \textbf{109}, no.6, 065207 (2024).

\bibitem{Weiss:2023yoj}
B.~Weiss, J.~F.~Paquet and S.~A.~Bass,
J. Phys. G \textbf{50}, no.6, 065104 (2023).

\bibitem{Ji:2022xzo}
Y.~Ji, H.~S.~Yuchi, D.~Soeder, J.~F.~Paquet, S.~A.~Bass, V.~R.~Joseph, C.~F.~J.~Wu and S.~Mak,
[arXiv:2209.13748 [stat.ME]].

\bibitem{Heffernan:2022swr}
M.~Heffernan, C.~Gale, S.~Jeon and J.~F.~Paquet,
Acta Phys. Polon. Supp. \textbf{16}, no.1, 1-A151 (2023).

\bibitem{Liyanage:2022byj}
D.~Liyanage, Y.~Ji, D.~Everett, M.~Heffernan, U.~Heinz, S.~Mak and J.~F.~Paquet,
Phys. Rev. C \textbf{105}, no.3, 034910 (2022).

\bibitem{JETSCAPE:2021rfx}
J.~Mulligan \textit{et al.} [JETSCAPE],
[arXiv:2106.11348 [nucl-th]].

\bibitem{JETSCAPE:2021ehl}
S.~Cao \textit{et al.} [JETSCAPE],
Phys. Rev. C \textbf{104}, no.2, 024905 (2021).

\bibitem{Bernhard:2019bmu}
J.~E.~Bernhard, J.~S.~Moreland and S.~A.~Bass,
Nature Phys. \textbf{15}, no.11, 1113-1117 (2019).



\bibitem{Moreland:2018gsh}
J.~S.~Moreland, J.~E.~Bernhard and S.~A.~Bass,
Phys. Rev. C \textbf{101}, no.2, 024911 (2020).

\bibitem{Moreland:2018jos}
J.~S.~Moreland, J.~E.~Bernhard and S.~A.~Bass,
Nucl. Phys. A \textbf{982}, 503-506 (2019).

\bibitem{Xu:2017obm}
Y.~Xu, J.~E.~Bernhard, S.~A.~Bass, M.~Nahrgang and S.~Cao,
Phys. Rev. C \textbf{97}, no.1, 014907 (2018).

\bibitem{Ke:2017xyi}
W.~Ke, J.~S.~Moreland, J.~E.~Bernhard and S.~A.~Bass,
Nucl. Part. Phys. Proc. \textbf{289-290}, 483-486 (2017).

\bibitem{Xu:2017nmc}
Y.~Xu, S.~Cao, M.~Nahrgang, J.~E.~Bernhard and S.~A.~Bass,
Nucl. Part. Phys. Proc. \textbf{289-290}, 257-260 (2017).

\bibitem{Auvinen:2017fjw}
J.~Auvinen, J.~E.~Bernhard, S.~A.~Bass and I.~Karpenko,
Phys. Rev. C \textbf{97}, no.4, 044905 (2018).

\bibitem{Xu:2017hgt}
Y.~Xu, M.~Nahrgang, J.~E.~Bernhard, S.~Cao and S.~A.~Bass,
Nucl. Phys. A \textbf{967}, 668-671 (2017).

\bibitem{Auvinen:2017pny}
J.~Auvinen, I.~Karpenko, J.~E.~Bernhard and S.~A.~Bass,
Nucl. Phys. A \textbf{967}, 784-787 (2017).

\bibitem{Bernhard:2017vql}
J.~E.~Bernhard, J.~S.~Moreland and S.~A.~Bass,
Nucl. Phys. A \textbf{967}, 293-296 (2017).

\bibitem{Ke:2016jrd}
W.~Ke, J.~S.~Moreland, J.~E.~Bernhard and S.~A.~Bass,
Phys. Rev. C \textbf{96}, no.4, 044912 (2017).

\bibitem{Auvinen:2016tgi}
J.~Auvinen, J.~E.~Bernhard and S.~A.~Bass,
Acta Phys. Polon. Supp. \textbf{10}, 455 (2017).

\bibitem{Bernhard:2015hxa}
J.~E.~Bernhard, P.~W.~Marcy, C.~E.~Coleman-Smith, S.~Huzurbazar, R.~L.~Wolpert and S.~A.~Bass,
Phys. Rev. C \textbf{91}, no.5, 054910 (2015).

\bibitem{Jahan:2024wpj}
S.~A.~Jahan, H.~Roch and C.~Shen,
Phys. Rev. C \textbf{110}, no.5, 054905 (2024).

\bibitem{JETSCAPE:2023nuf}
A.~Mankolli \textit{et al.} [JETSCAPE],
EPJ Web Conf. \textbf{296}, 05010 (2024).

\bibitem{Shen:2023pgb}
C.~Shen, B.~Schenke and W.~Zhao,
EPJ Web Conf. \textbf{296}, 14001 (2024).

\bibitem{Mantysaari:2022ypp}
H.~M\"antysaari, B.~Schenke, C.~Shen and W.~Zhao,
Acta Phys. Polon. Supp. \textbf{16}, no.1, 1-A33 (2023).

\bibitem{Mantysaari:2022ffw}
H.~M\"antysaari, B.~Schenke, C.~Shen and W.~Zhao,
Phys. Lett. B \textbf{833}, 137348 (2022).

\bibitem{JETSCAPE:2020mzn}
D.~Everett \textit{et al.} [JETSCAPE],
Phys. Rev. C \textbf{103}, no.5, 054904 (2021).

\bibitem{JETSCAPE:2020avt}
J.~F.~Paquet \textit{et al.} [JETSCAPE],
Nucl. Phys. A \textbf{1005}, 121749 (2021).

\bibitem{Pratt:2015zsa}
S.~Pratt, E.~Sangaline, P.~Sorensen and H.~Wang,
Phys. Rev. Lett. \textbf{114}, 202301 (2015).

\bibitem{Sangaline:2015isa}
E.~Sangaline and S.~Pratt,
Phys. Rev. C \textbf{93}, no.2, 024908 (2016).

\bibitem{Novak:2013bqa}
J.~Novak, K.~Novak, S.~Pratt, J.~Vredevoogd, C.~Coleman-Smith and R.~Wolpert,
Phys. Rev. C \textbf{89}, no.3, 034917 (2014).

\bibitem{Li:2025iqq}
B.~A.~Li and W.~J.~Xie,
[arXiv:2501.02579 [nucl-th]].


\bibitem{Alqahtani:2024ejg}
M.~Alqahtani, R.~S.~Bhalerao, G.~Giacalone, A.~Kirchner and J.~Y.~Ollitrault,
Phys. Rev. C \textbf{110}, no.6, 064906 (2024).

\bibitem{Wang:2024ahj}
J.~M.~Wang, X.~G.~Deng, W.~J.~Xie, B.~A.~Li and Y.~G.~Ma,
[arXiv:2406.07051 [nucl-th]].

\bibitem{Li:2023ydi}
B.~A.~Li and W.~J.~Xie,
Nucl. Phys. A \textbf{1039}, 122726 (2023).

\bibitem{Giannini:2022bkn}
A.~V.~Giannini \textit{et al.} [ExTrEMe],
Phys. Rev. C \textbf{107}, no.4, 044907 (2023).

\bibitem{Jia:2013tja}
J.~Jia and S.~Mohapatra,
Phys. Rev. C \textbf{88}, no.1, 014907 (2013).

\bibitem{Morfouace:2019jky}
P.~Morfouace et al., 
Phys. Lett. B \textbf{799}, 135045 (2019).

\bibitem{RasmussenGP}
C.E. Rasmussen and C.K.I. Williams, Gaussian Processes for Machine Learning, the MIT Press, 2006,
ISBN 026218253X.

\bibitem{WilliamsRasmussen96}
Williams, C.K.I. and Rasmussen, C.E., {\it Gaussian Processes for Regression} in {\it Advances in Neural Information Processing Systems} {\bf 8}, ed. by D.S. Touretzky, M.C. Mozer and M.E. Hasselmo, MIT Press (1996).

\bibitem{Mak24}
T. Tang, S. Mak and D. Dunson,
{\it Hierarchical Shrinkage Gaussian Processes: Applications to Computer Code Emulation and Dynamical System Recovery}
SIAM/ASA J. on Uncertainty Quantification {\bf 12}, no. 4, 10.1137 (2024).

\bibitem{Xiu2010}
D. Xiu, {\it Generalized Polynomial Chaos} in {\it Numerical Methods for Stochastic Computations}, Princeton University Press (2010)








\bibitem{lhs}
M.D. McKay, R.J. Beckman, R.J. and W.J. Conover.
Technometrics. 21 (2). American Statistical Association: 239–245.
(1979).
\end{thebibliography}
\end{document}